\documentclass[a4paper,twoside, 11pt]{article}
\usepackage{amssymb,amsfonts,amsmath}
\begin{document}
\newtheorem{defin}{~~~~Definition}
\newtheorem{prop}{~~~~Proposition}[section]
\newtheorem{remark}{~~~~Remark}[section]
\newtheorem{cor}{~~~~Corollary}
\newtheorem{theor}{~~~~Theorem}
\newtheorem{lemma}{~~~~Lemma}[section]
\newtheorem{ass}{~~~~Assumption}
\newtheorem{con}{~~~~Conjecture}
\newtheorem{concl}{~~~~Conclusion}
\renewcommand{\theequation}{\thesection.\arabic{equation}}
\newcommand {\trans} {^{\,\mid\!\!\!\cap}}
\newcommand{\vf}{\varphi}
\newcommand{\om}{\omega}
\newcommand{\bom}{\overline \omega}

\title{ Variational Approach to Differential Invariants of
Rank 2 Vector Distributions} 
\date{}
\author{Igor Zelenko\thanks{S.I.S.S.A., Via Beirut 2-4,
34014, Trieste, Italy; email: zelenko@sissa.it}} \maketitle 
\setcounter{equation}{0} 

\begin{abstract}
In the present paper we construct differential invariants 
for generic rank 2 vector distributions on $n$-dimensional 
manifold. In the case $n=5$ (the first case containing 
functional parameters) E. Cartan found in 1910 the 
covariant fourth-order tensor invariant for such 
distributions, using his "reduction-prolongation" procedure 
(see \cite{cartan}). After Cartan's work the following 
questions remained open: first the geometric reason for 
existence of Cartan's tensor was not clear; secondly it was 
not clear how to generalize this tensor to other classes of 
distributions; finally there were no explicit formulas for 
computation of Cartan's tensor. Our paper is the first in 
the series of papers, where we develop an alternative 
approach, which gives the answers to the questions 
mentioned above. It is based on the investigation of 
dynamics of the field of so-called abnormal extremals 
(singular curves) of rank $2$ distribution and on the 
general theory of unparametrized curves in the Lagrange 
Grassmannian, developed in \cite{jac1},\cite{jac2}. In this 
way we construct the {\it 
 fundamental form} and
 the {\it projective Ricci curvature} of rank 2 vector distributions for 
arbitrary $n\geq 5$.
 For $n=5$ we give an explicit method for computation of these
 invariants and demonstrate it on several examples. 
In the next paper \cite{zelcart} we show that in the case 
$n=5$ our fundamental form coincides 
 with Cartan's tensor.
 \end{abstract}

\vskip .2in


\section{Introduction}
\setcounter{equation}{0}
 \indent 

 Rank $l$ vector distribution $D$ on the $n$-dimensional
manifold $M$ or $(l,n)$-distribution (where $l<n$) is by 
definition a $l$-dimensional subbundle of the tangent 
bundle $TM$. In other words, for each point $q\in M$ a 
$l$-dimensional subspace $D(q)$ of the tangent space $T_qM$ 
is chosen and $D(q)$ depends smoothly on $q$. Two vector 
distributions $D_1$ and $D_2$ are called equivalent, if 
there exists a diffeomorphism $F:M\mapsto M$ such that 
$F_*D_1(q)=D_2(F(q))$ for any $q\in M$. Two germs of vector 
distributions $D_1$ and $D_2$ at the point $q_0\in M$ are 
called equivalent, if there exist neighborhoods $U$ and 
$\tilde U$ of $q_0$ and a diffeomorphism $F:U\mapsto \tilde 
U$ such that 
$$\begin{array}{c}F_*D_1(q)=D_2(F(q)),\,\,\forall q\in 
U;\\F(q_0)=q_0.\end{array}$$ Our goal is to construct 
invariants of distributions w.r.t. this equivalence 
relation in order to see if two germs of distributions are 
equivalent or not. Distributions are associated with 
Pfaffian systems and with control systems linear in the 
control. So invariants of distributions are also invariants 
of the corresponding Pfaffian systems and state-feedback 
invariants of the corresponding control systems.

An obvious (but very rough in the most cases) invariant of 
distribution $D$ at $q$ is so-called {\it small growth 
vector} at $q$: it is the tuple $$\bigl(\dim D(q),\dim 
D^2(q),\dim D^3(q),\ldots\bigr),$$ where $D^j$ is the 
$j$-th power of the distribution $D$, i.e., 
$D^j=D^{j-1}+[D,D^{j-1}]$.

Let us roughly estimate the "number of parameters" in the 
considered equivalence problem. The set of $l$-dimensional 
subspaces in $\mathbb{R}^n$ forms $l(n-l)$-dimensional 
manifold. Therefore, if the coordinates on $M$ are fixed 
then the rank $l$ distribution can be defined by $l(n-l)$ 
functions of $n$ variables. The group of the coordinate 
changes 
on $M$ is parameterized by $n$ functions of $n$ variables. 
So, by a coordinate change one can "normalize", in general, 
only $n$ functions among those $l(n-l)$ functions, defining 
the distribution. Thus, one may expect that the set of 
classes of equivalent germs of rank $l$ distributions 
can be parameterized by $l(n-l)-n$ arbitrary germs of 
functions of $n$ variables (see \cite{versh1} or survey 
\cite{versh2}, subsection 2.7, for precise statements). 
According to this, in the case $l=2$ the functional 
invariant should appear starting from $n=5$. It is well 
known that in the low dimensions $n=3$ or $4$ all generic 
germs of rank 2 distributions are equivalent. (Darboux's 
theorem in the case $n=3$, small growth vector $(2,3)$ and 
Engel's theorem in the case $n=4$, small growth vector 
$(2,3,4)$, see, for example, 
\cite{bryantbook}, \cite{zhi}).

The case of $(2,5)$-distributions with small growth vector 
$(2,3,5)$ was treated by E. Cartan in \cite{cartan} by 
ingenious use of his "reduction-prolongation" procedure. In 
particular, he constructed invariant homogeneous polynomial 
of degree $4$ on each plane $D(q)$ (we will call it 
Cartan's tensor). If the roots of the projectivization of 
this polynomial are different, then taking their 
cross-ratio one obtains functional invariant of the 
distribution $D$. 

After the mentioned work of E. Cartan the following 
questions remained open: first the geometric reason for 
existence of Cartan's tensor was not clear (the tensor was 
obtained by very sophisticated algebraic manipulations) and 
the true analogs of this tensor in Riemannian geometry were 
not found; secondly it was not clear how to generalize this 
tensor to other classes of distributions ; finally there 
were no explicit formulas for computation of Cartan's 
tensor 
(in order to compute this tensor for concrete distribution, 
one had to repeat Cartan's "reduction-prolongation" 
procedure for this distribution from the very beginning, 
which is rather difficult task). 

In the present paper we develop alternative, more geometric 
method for construction of functional invariants of generic 
germs of $(2,n)$-distribution for arbitrary $n\geq 5$, 
which allows to give the answers to the questions mentioned 
above. It is based on new, variational approach for 
constructing of differential invariants for 
geometric structures (feedback invariants for control 
systems) proposed recently by A. Agrachev (see 
\cite{agrgam1}, \cite{agrachev}, and also Introduction to 
\cite{jac1}). For rank 2 distributions this approach can be 
described as follows (the presentation here is closed to 
those given in Introduction to \cite{jac1}): 

First for $(2,n)$-distributions ($n\geq 4$) with small 
growth vector of the type $(2,3,4\, {\rm or}\, 5,\ldots)$ 
one can distinguish special (unparametrized) curves in the 
cotangent bundle $T^*M$ of $M$. For this let 
$\pi:T^*M\mapsto M$ be the canonical projection. Let 
$\sigma$ be standard symplectic structure on $T^*M$, 
namely, for any $\lambda\in T^*M$, $\lambda=(p,q)$, $q\in 
M$,$p\in T_q^*M$ let 
\begin{equation}
\label{stands}
\sigma(\lambda)(\cdot)=-d\, p\bigl(\pi_*(\cdot)\bigr)
\end{equation}
(here we prefer the sign ''-'' in the righthand side , 
although usually one defines the standard symplectic form 
on $T^*M$ without this sign). Denote by 
$(D^l)^{\perp}\subset T^*M$ the annihilator of the $l$th 
power $D^l$, namely 
\begin{equation}
\label{annihil} (D^l)^{\perp}= \{(q,p)\in T^*M:\,\, 
p\cdot v=0\,\,\forall v\in D^l(q)\}.
\end{equation}
The set $D^\perp$ is codimension $2$ submanifold of $T^*M$. 
Consider the restriction $\sigma|_{D^\perp}$ of the form 
$\sigma$ on $D^\perp$. It is not difficult to check that 
(see, for example \cite{zel}, section 2): the set of 
points, where the form $\sigma|_{D^\perp}$ is degenerated, 
coincides with $(D^2)^\perp$; the set 
$(D^2)^\perp\backslash(D^3)^\perp$ is codimension $1$ 
submanifold of $D^\perp$; for each $\lambda\in 
(D^2)^\perp\backslash(D^3)^\perp$ the kernel of 
$\sigma|_{D^\perp}(\lambda)$ is two-dimensional subspace of 
$T_\lambda D^\perp$, which is transversal to $T_\lambda 
(D^2)^\perp$. Hence $\forall\lambda\in 
(D^2)^\perp\backslash(D^3)^\perp$ we have $$ {\rm 
ker}\,\sigma|_{(D^2)^\perp}(\lambda)= {\rm 
ker}\,\sigma|_{D^\perp}(\lambda) \cap T_\lambda (D^2)^\perp 
$$ It implies that these kernels form line distribution in 
$(D^2)^\perp\backslash(D^3)^\perp$ and define a {\it 
characteristic 1-foliation} $Ab_D$ of 
$(D^2)^\perp\backslash(D^3)^\perp$. Leaves of this 
foliation will be called {\it characteristic curves} of 
distribution $D$. Actually these characteristic curves are 
so-called regular abnormal extremals of $D$ (see 
\cite{suss},\cite{zel} and also \cite {agrsar}, where such 
curves are called abnormal extremals, satisfying the strong 
generalized Legendre-Clebsh condition). 

\begin{remark}
\label{PMP} {\rm 
 The 
term abnormal extremal comes from Pontryagin Maximum 
Principle in Optimal Control. 
Defining on the set of all curves tangent to $D$ some 
functional (for example, length w.r.t. some Riemannian 
metric on $M$), one can consider the corresponding optimal 
control problem with fixed endpoints. Abnormal extremals 
are the extremals of this problem with vanishing Lagrange 
multiplier near the functional, so they do not depend on 
the functional but on the distribution $D$ itself. Projections 
of abnormal extremals to the base manifold $M$ will be 
called {\it abnormal trajectories}. Conversely, an abnormal 
extremal projected to the given abnormal trajectory will be 
called its {\it lift}. If some lift of the abnormal 
trajectory is regular abnormal extremal, then this abnormal 
trajectory will be called regular. Again from Pontryagin 
Maximum Principle it follows that the set of all lifts of 
given abnormal trajectory can be provided with the 
structure of linear space. The dimension of this space is 
called {\it corank of the abnormal trajectory}}.$\Box$ 
\end{remark}

Further, for a given segment $\gamma$ of characteristic 
curve one can construct a special (unparametrized) curve of 
Lagrangian subspaces, called Jacobi curve, in the 
appropriate symplectic space. For this for any $\lambda\in 
(D^2)^\perp$ denote by ${\cal J}(\lambda)$ the following 
subspace of $T_\lambda(D^2)^\perp$ 
\begin{equation}
\label{prejac} {\cal J}(\lambda)= \bigl(T_\lambda 
(T^*_{\pi(\lambda)}M)+ 
\ker\sigma|_{D^\perp}(\lambda)\bigr)\cap T_\lambda 
(D^2)^\perp.
\end{equation}
Here $T_\lambda (T^*_{\pi(\lambda)}M)$ is tangent to the 
fiber $T^*_{\pi(\lambda)}M$ at the point $\lambda$ (or 
vertical subspace of $T_\lambda(T^*M)$).
Actually ${\cal J}$ is rank $(n-1)$ distribution on the manifold 
$(D^2)^\perp\backslash (D^3)^\perp$. 

 Let $O_\gamma$ be 
a neighborhood of $\gamma$ in $(D^2)^\perp$ such that 
\begin{equation}
\label{Ndef} N=O_\gamma /({Ab_D}|_{O_\gamma})
\end{equation}
 is a 
well-defined smooth manifold. The quotient manifold $N$ is 
a symplectic manifold endowed with a symplectic structure 
$\bar\sigma$ induced by $\sigma |_{(D^2)^\perp}$. Let $\phi 
:O_\gamma\to N$ be the canonical projection on the factor. 
It is easy to check that $\phi_*\bigl({\cal 
J}(\lambda)\bigr)$ is a Lagrangian subspace of the 
symplectic space $T_{\gamma}N$, $\forall\lambda\in\gamma$. 
 Let $L(T_\gamma 
N)$ be the Lagrangian Grassmannian of the symplectic space 
$T_\gamma N$, i.e., $$L(T_\gamma N)=\{\Lambda\subset 
T_\gamma N:\Lambda^\angle=\Lambda\},$$ where 
$\Lambda^\angle$ is the skew-symmetric complement of the 
subspace $\Lambda$, $$\Lambda^\angle=\{v\in T_\gamma 
N:\bar\sigma(v,\Lambda)=0\}.$$ {\it Jacobi curve of the 
characteristic curve (regular abnormal extremal)} $\gamma$ 
is the mapping 
\begin{equation}
\label{jacurve} \lambda\mapsto 
J_\gamma(\lambda)\stackrel{def}{=} \phi_*\Bigl({\cal 
J}(\lambda) 
\Bigr), \quad\lambda\in\gamma , 
\end{equation}
from $\gamma$ to $L(T_\gamma N)$.   

\begin{remark}
\label{rigid} {\rm In \cite{agrgam1} and \cite{agrachev} 
Jacobi curves of extremals were constructed in purely 
variational way using the notion of Lagrangian derivative 
(${\cal L}$-derivative) of the endpoint map associated with 
geometric structure (control system). The reason to call 
these curves Jacobi curves is that they can be considered 
as generalization of spaces of "Jacobi fields" along 
Riemannian geodesics: in terms of these curves one can 
describe some optimality properties of corresponding 
extremals. Namely,   
if the Jacobi curve of the abnormal extremal is simple 
curve in Lagrange Grassmannian, 
 then the corresponding abnormal trajectory is 
$W_\infty^1$-isolated (rigid) curve in the space of all 
curves tangent to distribution $D$ with fixed endpoints ( 
the curve in Lagrange Grassmannian is called simple if one 
can choose Lagrangian subspace transversal to each Lagrange 
subspace belonging to the image of the curve). This result can be found in \cite{zel}. In different but equivalent form it is contained already in \cite{agrsar}. Moreover, if some Riemannian metric is given on $M$, then under the same conditions on  the Jacobi curve the corresponding abnormal trajectory is the shortest  among all curves tangent to distribution $D$, connecting 
its endpoints and sufficiently  closed to this abnormal trajectory in $W_1^1$-topology 
(see \cite {agrmin}) and even in $C$-topology (see \cite{agrc0}).}$\Box$ 
\end{remark}

Jacobi curves are invariants of the distribution $D$. They 
are unparametrized curves in the Lagrange Grassmannians.
 In 
\cite{jac1} for any curve of so-called constant {\it 
weight} in Lagrange Grassmannian we construct {\it the 
canonical projective structure} and the following two 
invariants w.r.t. the action of the linear Symplectic Group 
and reparametrization : a special degree 4 differential, 
{\it fundamental form}, and a special function, {\it 
projective Ricci curvature}. The next steps are to 
interpret the condition for Jacobi curve of regular 
abnormal extremal of distribution to be of constant weight 
in terms of distribution, to pass from the mentioned 
invariants defined on single Jacobi curve of each regular 
abnormal extremal of distribution to the corresponding 
invariants 
of distribution itself, and to investigate these 
invariants. These steps are the essence of the present 
work. 

The paper is organized as follows. In section 2 first we 
find under what assumption on germ of $(2,n)$-distribution 
at $q_0$ one can apply the general theory of unparametrized 
curves in the Lagrange Grassmannians to the Jacobi curve of 
its regular abnormal extremals. In few words this 
assumption can be described as follows: there is at least 
one germ of regular abnormal trajectory of distribution 
passing through $q_0$ and having corank 1 (see Remark 
\ref{PMP} for definitions). It is easy to see that the set 
of germs of $(2,n)$-distributions satisfying the last 
assumption is generic. In particular, in the case $n=5$ and 
$n=6$ the germs with the maximal possible small growth 
vector (namely, $(2,3,5)$ and $(2,3,5,6)$ respectively) 
satisfy this assumption. 

Further, for generic germ of $(2,n)$-distribution at $q_0$ 
we construct a fundamental form. By fundamental form at the 
point $q\in M$ we mean a special degree $4$ homogeneous 
rational function defined, up to multiplication on positive 
constant, on the linear space 
\begin{equation}
\label{D2q} (D^2)^\perp(q)=(D^2)^\perp\cap T^*_q M. 
\end{equation}
For germs of distributions at $q_0$ satisfying our assumption 
the set of points, where the fundamental 
form is defined, is open and dense in  some neighborhood of 
$q_0$. Later we show that for 
$(2,5)$-distribution with small growth vector $(2,3,5)$ the 
fundamental form at any point is polynomial, while for 
$n>5$ for generic $(2,n)$-distributions the fundamental 
form is a rational function, which is not a polynomial. 
Also, in the case of $(2,5)$-distribution with the small 
growth vector $(2,3,5)$ the fundamental form can be 
realized as degree 4 polynomial on the plane $D(q)$ for all 
$q\in M$ (we call it {\it tangential fundamental form}), 
i.e., it is an object of the same nature as Cartan's 
tensor. Further we describe the projective Ricci curvature 
of distribution, which is a function, defined on the subset 
of $(D^2)^\perp\backslash (D^3)^\perp$, where the 
fundamental form does not vanish. Note that the notion of 
projective Ricci curvature is new even for $n=5$. Using 
this notion, we construct, in addition to fundamental 
form, a special degree $10$ homogeneous rational function, 
defined, up to multiplication on positive constant, on 
$(D^2)^\perp(q)$ for any $q\in M$ (for $n=5$ this function 
is again polynomial). 

In section 3 we restrict ourselves to the case of 
$(2,5)$-distribution with small growth vector $(2,3,5)$. 
Using the notion of canonical moving frame of rank 1 curve 
in Lagrange Grassmannian introduced in \cite{jac1} (section 
6) and the structural equation for this frame derived in 
\cite{jac2} ( section 2) we obtain explicit formulas for 
fundamental form and projective curvature. It allows us to 
prove that in the considered case the fundamental form is a 
polynomial on each $(D^2)^{\perp}(q)$ (defined up to 
multiplication on positive constant). We apply the obtained 
formulas for several examples. In particular, we calculate 
our invariants for distribution generated by rolling of two 
spheres of radiuses $r$ and $\hat r$ ($r\leq \hat r$) 
without slipping and twisting. We show that the fundamental 
form of such distribution is equal to zero iff $\frac{\hat 
r}{r}=3$ and that the distributions with different ratios 
$\frac{\hat r} {r}$ are not equivalent. Also we give some 
sufficient conditions for rigidity of abnormal trajectories 
in terms of canonical projective structure and fundamental 
form on it.

 Finally, in section 4 we demonstrate that for $n>5$
 generically  the fundamental form 
 is not polynomial on the fibers $(D^2)^\perp(q)$. 
 It follows from the fact that fundamental form of the curve in 
 Lagrange Grassmannian
 always has singularities at the points of jump of the weight.

 In the next paper \cite{zelcart} we 
  prove that for $n=5$ our tangential fundamental form coincides 
  (up to constant factor $-35$) with Cartan's tensor. 
For this  we obtain another formula for fundamental form in 
 terms of structural functions of any frame naturally adapted to
  distribution  and  express 
  our fundamental form in terms of structural functions of special adapted 
  frame, distinguished by Cartan during the reduction. 
 
In the forthcoming paper \cite{zelprep} we construct the 
canonical frame of $(2,5)$- distribution on the subset of 
$(D^2)^{\perp}\backslash(D^3)^\perp$, where the fundamental 
form is not zero. It gives the way to check whether two 
different distribution are equivalent. Also we investigate 
distributions with constant projective Ricci curvature and 
big group of symmetries, giving models and proving 
uniqueness results ( we announce some of these results at 
the end of the section 3, see Theorems \ref{clastheor1}). 

Finally note that the approach of the present paper after 
some modifications can be applied to construction of 
invariant of other classes of distributions. In particular, 
the case of corank 2 distributions (i.e., when $n-l=2$) 
will be treated in the nearest future.

{\it Acknowledgments.} I would like to thank professor A. 
 Agrachev for his constant attention to this work, 
 stimulating discussions, and valuable advises.

\section{Fundamental 
form and projective Ricci curvature of rank 2 distribution } \indent 
\setcounter{equation}{0}

{\bf 2.1 Preliminary} Let $W$ be $2m$-dimensional linear 
space and $G_m(W)$ be the set of all $m$-dimensional 
subspaces of $W$ (i.e., the Grassmannian of 
half-dimensional subspaces). Below we give definitions of 
weight and rank of the curve in $G_m(W)$ and describe 
briefly the construction of fundamental form and projective 
Ricci curvature for the curve of constant weight in 
$G_m(W)$ (for the details see \cite{jac1}), which are invariants w.r.t. the action of General Linear Group $GL(W)$. Since any curve 
of Lagrange subspaces w.r.t. some symplectic form on $W$ is 
obviously the curve in $G_m(W)$, all constructions below 
are valid for the curves in Lagrange Grassmannian. 

For given $\Lambda\in G_m(W)$ denote by 
$\Lambda^\pitchfork$ the set of all $m$-dimensional 
subspaces of $W$ transversal to $\Lambda$, 
$$\Lambda^\pitchfork=\{\Gamma\in 
G_m(W):W=\Gamma\oplus\Lambda\}=\{\Gamma\in 
G_m(W):\Gamma\cap\Lambda=0\}$$ Fix some $\Delta\in 
\Lambda^\pitchfork$. Then for any subspace $\Gamma\in 
\Lambda^\pitchfork$ there exist unique linear mapping from 
$\Delta$ to $\Lambda$ with graph $\Gamma$. We denote this 
mapping by $\langle\Delta,\Gamma,\Lambda\rangle$. So, 
$$\Gamma=\{v+\langle\Delta,\Gamma,\Lambda\rangle 
v|v\in\Delta\}.$$ Choosing the bases in $\Delta$ and 
$\Lambda$ one can assign to any $\Gamma\in 
\Lambda^\pitchfork$ the matrix of the mapping
$\langle\Delta,\Gamma,\Lambda\rangle$ w.r.t these bases. In 
this way we define the coordinates on the set 
$\Lambda^\pitchfork$. 

\begin{remark}
\label{lagr} {\rm Assume that $W$ is endowed with some 
symplectic form $\bar\sigma$ and $\Delta, \Lambda$ are 
Lagrange subspaces w.r.t. $\bar\sigma$. Then the map 
$v\mapsto \bar\sigma(v, \cdot)$, $v\in \Delta$, defines the 
canonical isomorphism between $\Delta$ and $\Lambda^*$. It 
is easy to see that $\Gamma$ is Lagrange subspace iff the 
mapping $\langle\Delta,\Gamma,\Lambda\rangle$, considered 
as the mapping from $\Lambda^*$ to $\Lambda$, is 
self-adjoint.} $\Box$ \end{remark} 

Let $\Lambda(t)$ be a smooth curve in $G_m(W)$ defined on 
some interval $I\subset\mathbb{R}$. We are looking for 
invariants of $\Lambda(t)$ by the action of $GL(W)$. 
We say that the curve $\Lambda(\cdot)$ is {\it ample at} 
$\tau$ if $\exists s>0$ such that for any representative 
$\Lambda^s_\tau(\cdot)$ of the $s$-jet of $\Lambda(\cdot)$ 
at $\tau$, $\exists t$ such that 
$\Lambda^s_\tau(t)\cap\Lambda(\tau)=0$. The curve 
$\Lambda(\cdot)$ is called {\it ample} if it is ample at 
any point. This is an intrinsic definition of an ample 
curve. In coordinates this definition takes the following 
form: if in some coordinates the curve $\Lambda(\cdot)$ is 
a curve of matrices $t\mapsto S_t$, then $\Lambda(\cdot)$ 
is ample at $\tau$ if and only if the function $t\mapsto 
det(S_t-S_\tau)$ has a root of {\sl finite order} at 
$\tau$. 

\begin{defin}
The order of zero of the function $t\mapsto 
det(S_t-S_\tau)$ at $\tau$, where $S_t$ is a coordinate 
representation of the curve $\Lambda(\cdot)$, is called a 
weight of the curve $\Lambda(\cdot)$ at $\tau$. 
\end{defin}  
It is clear that the weight of $\Lambda(\tau)$ is integral 
valued upper semicontinuous functions of $\tau$. Therefore 
it is locally constant on the open dense subset of the 
interval of definition $I$. 

Now suppose that the curve has the constant weight $k$ on 
some subinterval $I_1\subset I$. It implies that for all 
two parameters $t_0$,$t_1$ in $I_1$ sufficiently such that 
$t_0\neq t_1$ , one has 
 $$
 \Lambda(t_0)\cap\Lambda(t_1)=0.
 $$
 Hence for such $t_0$, $t_1$ the following linear mappings
 \begin{eqnarray} ~&~
\frac{d}{d s}\langle\Lambda(t_0),\Lambda(s),\Lambda(t_1)
 \rangle\Bigl|_{s=t_0}\Bigr .:\Lambda(t_0)\mapsto\Lambda(t_1),\label{op1}\\ 
 ~&~\frac{d}{d s}\langle\Lambda(t_1),\Lambda(s),\Lambda(t_0)
 \rangle\Bigl|_{s=t_1}\Bigr 
 .:\Lambda(t_1)\mapsto\Lambda(t_0)\label{op2}
 \end{eqnarray}
 are well defined. Taking composition of mapping (\ref{op2}) 
 with mapping  (\ref{op1}) we obtain the operator 
 from the subspace $\Lambda(t_0)$ to itself, which is 
 actually the infinitesimal cross-ratio of two points
 $\Lambda(t_i)$, $i=0,1$, together with two tangent vectors 
 $\dot\Lambda(t_i)$, $i=0,1$, at these points in $G_m(W)$ 
 (see \cite{jac1} for the details).

 \begin{theor} (see \cite{jac1}, Lemma 4.2) If the curve has the constant weight $k$ on 
some subinterval $I_1\subset I$, then the following 
asymptotic holds 
 \begin{eqnarray}
  ~&~{\rm tr}\Bigl( \frac{d}{ds}\langle\Lambda(t_1),\Lambda(s),\Lambda(t_0)
 \rangle\Bigl|_{s=t_1}\Bigr.\circ
  \frac{d}{ds}\langle\Lambda(t_0),
 \Lambda(s),\Lambda(t_1)
 \rangle\Bigl|_{s=t_0}\Bigr.\Bigr)=\nonumber\\
 ~~\label{mainass}\\
 ~&~-\frac {k}{(t_0-t_1)^2}-g_\Lambda(t_0,t_1),\nonumber
 \end{eqnarray}
 where $g_{_{\Lambda}}(t,\tau)$ 
  is a smooth function in the neighborhood of diagonal 
  $\{(t,t)| t\in I_1\}$. 
\end{theor}

The function $g_{_{\Lambda}}(t,\tau)$ is "generating" 
function for invariants of the parametrized curve by the 
action of $GL(2m)$. The first coming invariant of the 
parametrized curve, {\it the generalized Ricci curvature}, is 
just $g_{_{\Lambda}}(t,t)$, 
 the value of $g_{_{\Lambda}}$  at the diagonal. 

In order to obtain invariants for unparametrized curves 
(i.e., for one-dimensional submanifold of $G_m(W)$) we use 
a simple reparametrization rule for a function 
$g_{_{\Lambda}}$. Indeed, let 
$\vf:\mathbb{R}\mapsto\mathbb{R}$ be a smooth monotonic 
function. It follows directly from (\ref{mainass}) that 
\begin{equation}\label{chainrul} g_{_{\Lambda\circ 
\vf}}(t_0,t_1)=\dot\vf(t_0)\dot\vf(t_1)g_{_{\Lambda}} 
(\vf(t_0),\vf(t_1))+k\left(\frac{\dot\vf(t_0)\dot\vf(t_1)} 
{(\vf(t_0)-\vf(t_1))^2}-\frac{1}{(t_0-t_1)^2}\right). 
\end{equation}
 In particular, putting $t_0=t_1=t$, one obtains the following reparametrization rule for the generalized Ricci curvature 
\begin{equation}
 \label{rhorep}
 g_{_{\Lambda\circ 
\vf}}(t,t)=\dot\vf(t)^2g_{_{\Lambda}}(\vf(t),\vf(t))+ 
\frac{k}{3}\mathbb{S}(\vf),
\end{equation}
 where $\mathbb{S}(\vf)$ is a 
Schwarzian derivative of $\vf$, \begin{equation} 
\label{sch} \mathbb S(\varphi)= 
\frac{1}{2}\frac{\varphi^{(3)}}{\varphi'}- \frac{3}{4}\Bigl 
(\frac{\varphi''}{\varphi'}\Bigr)^2=\frac 
{d}{dt}\Bigl(\frac {\varphi''}{2\,\varphi'}\Bigr) 
-\Bigl(\frac{\varphi''}{2\,\varphi'}\Bigr)^2.
\end{equation}

From (\ref{rhorep}) it follows that the class of local 
parametrizations, in which the generalized Ricci curvature 
is identically equal to zero, defines a {\it canonical 
projective structure} on the curve (i.e., any two 
parametrizations from this class are transformed one to 
another by M{\"obius transformation). This parametrizations 
are called projective. From (\ref{chainrul}) it follows 
that if $t$ and $\tau$ are two projective parametrizations 
on the curve $\Lambda(\cdot)$, 
$\tau=\vf(t)=\frac{at+b}{ct+d}$, and $g_\Lambda$ is 
generating function of $\Lambda(\cdot)$ w.r.t. parameter 
$\tau$ then 
\begin{equation}
\label{prefund} {\partial^2 
g_{_{\Lambda\circ\vf}}\over\partial 
t_1^2}(t_0,t_1)\Bigr|_{t_0=t_1=t}\Bigl.={\partial^2 
g_{_{\Lambda}}\over\partial 
\tau_1^2}(\tau_0,\tau_1)\Bigr|_{\tau_0=\tau_1=\vf(t)}\Bigl. 
(\vf'(t))^4, 
\end{equation}
which implies that the following degree four differential 
$$ {\cal A}=\frac{1}{2} {\partial^2 
g_{_{\Lambda}}\over\partial 
\tau_1^2}(\tau_0,\tau_1)\Bigr|_{\tau_0=\tau_1=\tau}\Bigl. 
(d\tau)^4$$ on the curve $\Lambda(\cdot)$ does not depend 
on the choice of the projective parametrization (by degree 
four differential on the curve we mean the following: for 
any point of the curve a degree $4$ homogeneous function is 
given on the tangent line to the curve at this point). This 
degree four differential is called a {\it fundamental form} 
of the curve. 

If $t$ is an arbitrary (not necessarily projective) 
parametrization on the curve $\Lambda(\cdot)$, then the 
fundamental form in this parametrization has to be of the 
form $A(t)(dt)^4$, where $A(t)$ is a smooth function, 
called the {\it density} of the fundamental form w.r.t. 
parametrization $t$. The density $A(t)$ can be expressed by 
generating function $g_\Lambda$ in the following way (see 
\cite{jac1}, Lemma 5.1): 

\begin{equation}
\label{densg} 
A(t)=
\frac{1}{2}{\partial^2\over\partial t_0^2 
}g_\Lambda(t_0,t_1)\Bigl|_{t_0=t_1=t}\Bigr.-\frac{3}{5k}g_\Lambda(t,t)^2-\frac{3}{20}\frac{d^2}{d 
t^2}g_\Lambda(t,t). 
\end{equation}

\begin{remark}
\label{ratrem} {\rm If $t\mapsto S_t$ is a coordinate 
representation of the curve $\Lambda(\cdot)$ in some 
parametrization $t$, then 
\begin{equation}
\label{gdet}
g_\Lambda(t_0,t_1)={\partial^2\over\partial t_0\partial t_1}
\ln\left(\frac{\det(S_{t_0}-S_{t_1})}{(t_0-t_1)^k}\right)
\end{equation}
(the proof of the last formula follows from \cite{jac1}, see relations (4.9),(4.11), and Lemma 4.2 there).
From this
and 
(\ref{densg}) it follows that for any $t_0$ the density $ 
A(t_0)$ w.r.t. parametrization $t$ of the fundamental form 
of $\Lambda(\cdot)$ at $t_0$ is equal to a rational 
expression w.r.t. some entries of the matrices $\dot 
S(t_0),\ddot S(t_0),\ldots,S^{(j)}(t_0)$ for some $j>0$.} 
$\Box$ 
\end{remark}

If the fundamental form ${\cal A}$ of the curve 
$\Lambda(\cdot)$ is not equal to zero at any point of 
$\Lambda(\cdot)$, then the canonical length element 
$\sqrt[4]{|{\cal A}|}$ is defined on $\Lambda(\cdot)$. The 
length with respect to this length element gives canonical, 
up to the shift, parametrization of the unparametrized 
curve. The Ricci curvature $\rho_n(\tau)$ w.r.t. this 
parametrization is a functional invariant of the 
unparametrized curve, 
which is called its {\it projective Ricci curvature}. 
\begin{remark}
\label{denscan} {\rm By construction, the density of the 
fundamental form in the canonical parameter is identically 
equal to 1.} $\Box$
\end{remark}
\begin{remark}
\label{m=1A=0} {\rm Note that in the case $m=1$ the 
fundamental form is always identically equal to zero (see 
\cite {jac1}, Lemma 5.2 there): in this case the only 
invariant of unparametrized curve in the corresponding 
Lagrange Grassmannian is the canonical projective structure 
on it.} $\Box$ 
\end{remark}

 Another important characteristic of the 
curve $\Lambda(\cdot)$ in $G_m(W)$ 
 is the rank of its velocities (or simply 
rank). In order to define it let
\begin{equation}
\label{primext} {\cal 
D}^{(1)}\Lambda(\tau)\stackrel{def}{=}\Lambda (\tau)+ 
\left\{v\in W: \begin{array}{l}\exists\,\, {\rm a}\,\,{\rm 
curve}\,\, 
 l(\cdot)\,\,{\rm in}\,\,  W\,\,{\rm such}\,\,{\rm that}\,\,\\
  l(t)\in \Lambda(t)\,\,\forall t 
\,\,{\rm  and}\,\,  v=\frac{d}{dt }l(t)|_{t=\tau}\end{array}\right \}
\end{equation} 
Then the rank $r(\tau)$ of $\Lambda(\cdot)$ at $\tau$ is defined as follows
\begin{equation}
\label{rank} r(\tau)\stackrel{def}{=}{\rm dim}\,{\cal 
D}^{(1)}\Lambda(\tau)- {\rm dim}\, \Lambda(\tau). 
\end{equation}
\begin{remark}
{\rm 
Note that the 
tangent space $T_\Lambda G_m(W)$ to any subspace $\Lambda\in G_m(W)$ 
can be identified with the space ${\rm Hom}\, (\Lambda,W/\Lambda)$ of 
linear mappings from $\Lambda$ to $W/\Lambda$. 
Namely, take a curve 
$\Lambda(t)\in G_k(W)$ with $\Lambda(0)=\Lambda$. Given some vector 
$l\in\Lambda$, take a curve $l(\cdot)$ in $W$ such that 
 $l(t)\in \Lambda (t)$ for all sufficiently small $t$ 
 and $l(0)=l$. Denote by ${\rm pr}:W\mapsto W/\Lambda$ the canonical 
 projection on the factor. 
 It is easy to see that the mapping $l\mapsto {\rm  pr}\,  
 l'(0)$ from $\Lambda$ to $W/\Lambda$ is linear mapping depending only on 
$\frac{d}{dt}\Lambda(0)$. In this way we identify  
$\frac{d}{dt}\Lambda(0)\in T_\Lambda G_k(W)$ with 
some element of ${\rm Hom}\, (\Lambda,W/\Lambda)$ 
(a  simple counting of 
dimension shows that these correspondence between $T_\Lambda G_k(W)$ and
 ${\rm Hom}\, (\Lambda,W/\Lambda)$ is a bijection). 
By construction, the rank of the curve $\Lambda(t)$ at the 
point $\tau$ in $G_m(W)$ is actually equal to the rank of 
the linear mapping corresponding to its velocity 
$\frac{d}{dt}\Lambda(t)$ at $\tau$.}$\Box$ 
\end{remark}

\begin{remark}
\label{tanglw} {\rm If $W$ is endowed with some symplectic 
form $\bar \sigma$ and $L(W)$ is the corresponding Lagrange 
Grassmannian, then the tangent space $T_\Lambda L(W)$ to 
any $\Lambda\in L(W)$ can be identified with the space of 
quadratic forms $Q(\Lambda)$ on the linear space $\Lambda$. 
Namely, let $\Lambda(t)$ and $l(t)$ be as in the previous 
remark (where $G_m(W)$ is substituted by $L(W)$). It is 
easy to see that the quadratic form $l\mapsto 
\bar\sigma\bigl(l'(0),l\bigr)$ depends only on 
$\frac{d}{dt}\Lambda(0)$. In this way we identify 
$\frac{d}{dt}\Lambda(0)\in T_\Lambda G_k(W)$ with some 
element of $Q(\Lambda)$ (a simple counting of dimension 
shows that these correspondence between $T_\Lambda L(W)$ 
and 
 $Q(\Lambda)$ is a bijection).} $\Box$
\end{remark}

 Using the identification in the previous remark one can define the notion of monotone curves in the Lagrange Grassmannian: the curve $\Lambda(t)$ in $L(W)$ is 
called
{\it nondecreasing (nonincreasing)} if its velocities $\frac{d}{dt}\Lambda(t)$ at any point are nonnegative (nonpositive) definite quadratic forms.

As we will see in the next subsection the rank of Jacobi 
curves of characteristic curves of rank $2$ distribution is 
identically equal to 1. There is a simple criterion for 
rank 1 curves in Lagrange Grassmannian to be of constant 
weight. To formulate it let us introduce inductively the 
following subspaces ${\cal D}^{(i)}\Lambda(\tau)$ in 
addition to ${\cal D}^{(1)}\Lambda(\tau)$: 

\begin{equation}
\label{primeiext} {\cal D}^{(i)}\Lambda(\tau)\stackrel{def} 
{=}{\cal D}^{(i-1)}\Lambda (\tau)+ \left\{v\in W: 
\begin{array}{l}\exists\,\, {\rm a}\,\,{\rm curve}\,\, 
 l(\cdot)\,\,{\rm in}\,\,  W\,\,{\rm such}\,\,{\rm that}\,\,\\
  l(t)\in {\cal D}^{(i-1)}\Lambda(t)\,\,\forall t 
\,\,{\rm  and}\,\,  v=\frac{d}{dt }l(t)|_{t=\tau}\end{array}\right \}
\end{equation} 
(we set $Der^{(0)}\Lambda(t)=\Lambda(t)$)

\begin{prop}
\label{constrank} The curve $\Lambda(\cdot)$ of constant 
rank 1 in Lagrange Grassmannian $L(W)$ of symplectic space 
$W$, ${\rm dim}\,\,W=2m$, has the constant finite weight in 
a neighborhood of the point $\tau$ iff 
\begin{equation}
\label{condconst} {\rm dim}\, {\cal 
D}^{(m)}\Lambda(\tau)=2m 
\end{equation}
In this case the weight is equal to $m^2$.
\end{prop}

The proof of the proposition can be easily obtained by 
application of some formulas and statements of section 6 
and 7 of \cite{jac1} (for example, formulas (6.15), (6.16), 
(6.18), (6.19), Proposition 4, and  Corollary 2  there).

Note also that from the fact that the rank of the curve is equal to $1$ it follows easily that
\begin{equation}
\label{dimjump} {\rm dim}\, {\cal D}^{(i)}\Lambda(\tau)- 
{\rm dim}\,{\cal D}^{(i-1)} \Lambda(\tau)\leq 1 
\end{equation}
Therefore the condition (\ref{condconst}) is equivalent to 
\begin{equation}
\label{condconst1} {\rm dim}\, {\cal 
D}^{(i)}\Lambda(\tau)=m+i, \,\, \forall i=1,\ldots m 
\end{equation}

{\bf 2.2 Properties of Jacobi curves of regular abnormal 
extremals of rank 2 distributions} In this subsection we 
find under what assumption on germ of $(2,n)$-distribution 
($n\geq 4$) with small growth vector of the type 
$(2,3,4\,{\rm or}\, 5,\ldots )$ one can apply the theory of 
previous subsection. 
 First note that Jacobi curve $J_\gamma$ of characteristic curve $\gamma$ 
 of distribution $D$
 defined by (\ref{jacurve}) is not 
ample, because all subspaces $J_\gamma(\lambda)$ have a 
common line. Indeed, let $\delta_a:T^*M\mapsto T^*M$ be the 
homothety by $a\neq 0$ in the fibers, namely, 

\begin{equation}
\label{homoth}
 \delta_a(p,q)=(ap,q), \quad q\in M,\,\,p\in 
T^*M.
\end{equation}
 Denote by $\vec e(\lambda)$ the following vector 
field called Euler field \begin{equation} \label{Euler} 
\vec e(\lambda)= \frac{\partial}{\partial a}\delta
_a(\lambda)\Bigl|_{a=1}\Bigr. 
\end{equation} 
\begin{remark}
\label{homrem}{\rm Obviously, if $\gamma$ is characteristic 
curve of $D$, 
then also $\delta_a(\gamma)$ is.} $\Box$ 
\end{remark}

It implies that the vectors $\phi_*(\vec 
e(\lambda))$ coincide for all $\lambda\in \gamma$, so the 
line 
\begin{equation}
\label{ega} E_\gamma\stackrel{def}{=}\{\mathbb R 
\phi_*(\vec e(\lambda))\}
\end{equation}
 is common for all subspaces 
$J_\gamma(\lambda)$, $\lambda\in \gamma$ (here, as in 
Introduction, $\phi 
:O_\gamma\to N$ is the canonical projection on the factor 
$N=O_\gamma /({Ab_D}|_{O_\gamma})$, where $O_\gamma$ is 
sufficiently small tubular neighborhood of $\gamma$ in 
$(D^2)^\perp$). 

Therefore it is natural to make an appropriate 
factorization by this common line $E_\gamma$. Namely, by 
above all subspaces $J_\gamma(\lambda)$ belong to 
skew-symmetric complement $E_\gamma^\angle$ of $E_\gamma$ 
in $T_\gamma N$. Denote by $p:T_\gamma N\mapsto T_\gamma 
N/E_\gamma$ the canonical projection on the factor-space. 
The mapping 
\begin{equation}
\label{redjac} \lambda\mapsto \widetilde 
J_\gamma(\lambda)\stackrel{def}{=}p(J_\gamma(\lambda)), 
\,\,\lambda\in\gamma 
\end{equation}
from $\gamma$ to $L(E_\gamma^\angle/E_\gamma)$ is called 
{\it reduced Jacobi curve} of characteristic curve 
$\gamma$. 
Note that 
\begin{equation}
\label{dimjac}
{\rm dim }\, \widetilde 
J_\gamma(\lambda)=n-3
\end{equation}

Now the question is at which points 
$\lambda\in \gamma$ the germ of reduced Jacobi curve has 
constant weight? The answer on this question can be easily done in terms of 
rank $(n-1)$ distribution  ${\cal J}$  defined by (\ref{prejac})  on 
$(D^2)^\perp\backslash (D^3)^\perp$.

First note that for any $\lambda\in\gamma$ one can make the 
following identification 
\begin{equation}
\label{tanid} 
T_\gamma N\sim  T_\lambda(D^2)^\perp/{\rm ker}\, \sigma|_{(D^2)^\perp}(\lambda).
\end{equation}

Take on $O_\gamma$ any vector field $H$ 
 tangent to characteristic $1$-foliation $Ab_D$ and without stationary points,
i.e.,  $H(\lambda)\in {\rm ker} \sigma|_{(D^2)^\perp}(\lambda)$, $H(\lambda)\neq 0$ for all $\lambda\in O_\gamma$. Then it is not hard to see that
under identification (\ref{tanid}) one has
\begin{equation}
\label{ident} \widetilde J_\gamma(e^{t H}\lambda)=(e^{-t H})_
*\Bigl( {\cal J}(e^{t 
H}\lambda) \Bigr) /{\rm span}( {\rm 
ker}\,\sigma|_{(D^2)^\perp}(\lambda),\vec e(\lambda) ) 
\end{equation} 
where $e^{t H}$ is the flow generated by the vector 
field $ H$.
Recall that for any vector field $\ell$ in 
$(D^2)^\perp$ one has 
\begin{equation}
\label{addt} \frac{d}{dt}\Bigl((e^{-t H})_*\ell 
\Bigr)
=(e^{-tH})_*[H,\ell] \end{equation}
Set ${\cal J}^{(0)}={\cal J}$
and define inductively
\begin{equation}
{\cal J}^{(i)}(\lambda)={\cal J}^{(i-1)}(\lambda)+
\{[H, V](\lambda): H\in {\rm ker}\,\sigma|_{(D^2)^\perp},
V\in {\cal J}^{(i-1)}\,\,{\rm are}\,\, 
{\rm vector}\,\, {\rm fields}\}
\label{Ji}
\end{equation} 
or shortly ${\cal J}^{(i)}={\cal J}^{(i-1)}+[{\rm 
ker}\,\sigma|_{(D^2)^\perp},{\cal J}^{(i-1)}]$. Then by 
definition of operation ${\cal D}^{(i)}$ (see 
(\ref{primeiext})) and formulas (\ref{ident}), (\ref{addt}) 
it follows that 

\begin{equation}
\label{DerJ} {\cal D}^{(i)}\widetilde 
J_\gamma(e^{tH}\lambda)|_{t=0}= {\cal 
J}^{(i)}(\lambda)/{\rm span}( {\rm 
ker}\,\sigma|_{(D^2)^\perp}(\lambda),\vec e(\lambda) ), 
\end{equation}
which in turn implies that
\begin{equation}
\label{derjJ} {\rm dim}\,{\cal D}^{(i)}\widetilde 
J_\gamma(e^{tH}\lambda)|_{t=0}- {\rm dim}\,{\cal D
}^{(i-1)}\widetilde J_\gamma(e^{tH}\lambda)|_{t=0}= {\rm 
dim}\,{\cal J}^{(i)}(\lambda)-{\rm dim}\, {\cal 
J}^{(i-1)}(\lambda). 
\end{equation}
\begin{prop}
\label{ranklem} 
 The (reduced) Jacobi curve of characteristic curve of $(2,n)$-distribution 
($n\geq 4$) with small growth vector of the type 
$(2,3,4\,{\rm or}\, 5,\ldots )$ is of rank 1 at any point 
and nondecreasing. 
 \end{prop}

{\bf Proof.} First show that the (reduced) Jacobi curve has 
rank 1 at any point. For this, according to (\ref{derjJ}), 
it is sufficient to prove that 
\begin{equation}
\label{equivrk}
{\rm dim}\,{\cal J}^{(1)}(\lambda)-{\rm dim}\, {\cal J}(\lambda) =1
\end{equation}

Let  $X_1$, $X_2$ be two vector fields, constituting
the basis of distribution $D$,
 i.e.,
\begin{equation}
\label{X12} D( q)={\rm span}(X_1(q), X_2(q))\,\,\,\,\forall 
q\in 
M.
\end{equation}
Since our study is local, we can always suppose that such 
basis exists, restricting ourselves, if necessary, on some 
coordinate neighborhood instead of whole $M$. 
Given the basis  $X_1$, $X_2$  one can construct a special 
vector field tangent to characteristic $1$-foliation 
$Ab_D$.
For this suppose that 
\begin{eqnarray}
&~&X_3=[X_1,X_2]\quad {\rm mod}\, D,\,\,\, 
X_4=[X_1,[X_1,X_2]]=[X_1,X_3]\quad {\rm mod}\, 
D^2,\nonumber 
\\ &~&~\label{x345} \\ 
&~&X_5=[X_2,[X_1,X_2]]=[X_2,X_3]\quad{\rm mod}\, D^2 
\nonumber 
\end{eqnarray}
Let us introduce ``quasi-impulses'' $u_i:T^*M\mapsto\mathbb 
R$, $1\leq i\leq 5$, 
\begin{equation}
\label{quasi25} u_i(\lambda)=p\cdot X_i(q),\,\,\lambda=(p,q),\,\,
q\in M,\,\, p\in T_q^* M
\end{equation}
For given function $G:T^*M\mapsto \mathbb R$ denote by 
$\vec G$ the corresponding Hamiltonian vector field defined 
by the relation $\sigma(\vec G,\cdot)=d\,G(\cdot)$.
Then it is easy to show (see, for example \cite{zel})
that
\begin{equation}
\label{ker25} \ker\sigma\Bigr|_{D^\perp}\Bigl.(\lambda)= 
{\rm span}(\vec
u_1(\lambda),\vec u_2(\lambda)),\quad \forall \lambda\in D^\perp, 
\end{equation} 
\begin{equation}
\label{foli25}
\ker\sigma\Bigr|_{(D^2)^\perp}\Bigl.(\lambda)= \mathbb{R}
\Bigl((u_4 \vec{u}_2-u_5\vec{u}_1)(\lambda)\Bigr),\quad 
\forall \lambda\in (D^2)^\perp\backslash (D^3)^\perp
\end{equation}
The last relation implies that 
the following vector field 
\begin{equation}
\label{ham25} \vec h_{_{X_1,X_2}}=u_4 \vec u_2- u_5 \vec 
u_1 
\end{equation}
is tangent to the characteristic $1$-foliation (this field 
is actually the restriction on $(D^2)^\perp$ of the 
Hamiltonian vector field of the function $h_{_{X_1,X_2}}= 
u_4 u_2-u_5 u_1$). 

Suppose that $\dim D^3(q)=5$ for any $q$ (the case, when 
$\dim D^3(q)=4$ for some $q$ can be treated similarly and 
it is left to the reader). Let us complete tuple 
$(X_1,X_2,X_3,X_4,X_5)$ to the local frame $X_1,\ldots, 
X_n$ on $M$. 
Similar to (\ref{quasi25})  define "quasi-impulses" 
$u_i:T^* M\mapsto 
\mathbb{R}$, $5< i\leq n$.

Denote 
\begin{eqnarray}
&~&
\partial_\theta=u_4\partial_{u_5}-u_5\partial_{u_4},\label{theta25}\\
 &~&{\cal X}=u_5\vec u_2+u_4\vec
 u_1-(u_4^2+u_5^2)\partial_{u_3},\label{X25}\\
 &~&F=\vec u_3+u_4\partial_{u_1}+u_5\partial_{u_2}\label{F25}
\end{eqnarray}
On $(D^2)^\perp \backslash (D^3)^\perp$ one has 
\begin{equation}
{\cal J}={\rm span} (\vec h_{_{X_1,X_2}},\vec e, {\cal 
X},\partial_\theta,\partial_{u_6}, 
\ldots,\partial_{u_n})\label{plambda} 
\end{equation}
 By direct computations, one can obtain
that
\begin{eqnarray}
&~& [\vec h_{_{X_1,X_2}},\partial_{u_i}]\in {\rm span}(\vec 
e, 
\partial_\theta,\partial_{u_6}, \ldots,\partial_{u_n}),\,\, 6\leq
i\leq n \label{part6}\\ &~&[\vec 
h_{_{X_1,X_2}},\partial_\theta]\equiv {\cal X} 
\,\,\,\Bigl({\rm mod}\bigl({\rm span}(\vec 
h_{_{X_1,X_2}},\vec e,\partial_\theta,\partial_{u_6}, 
\ldots,\partial_{u_n})\bigr)\Bigr),\label{pattheta}\\ 
&~& [\vec h_{_{X_1,X_2}}, {\cal X}]\equiv 
-(u_4^2+u_5^2)F\,\,\, ({\rm mod}\,\, {\cal J})
\label{hX}. 
\end{eqnarray}
From this and definition of ${\cal J}^{(1)}$ it follows that
\begin{equation}
\label{j1}
{\cal J}^{(1)}=\mathbb {R} F\oplus {\cal J},
\end{equation}
which implies (\ref{equivrk}). 

Finally, from (\ref{X25}), (\ref{F25}), and (\ref {hX}), it follows easily that  
$$\bar\sigma([h,{\cal X}],{\cal X})=(u_4^2+u_5^2)^2>0,$$
which implies that the (reduced) Jacobi curve is nondecreasing (see Remark \ref {tanglw} and the sentence after it).
$\Box$
\medskip

Proposition \ref{constrank}, relation (\ref{dimjump}), 
Proposition \ref{ranklem} and relation (\ref {derjJ}) imply 
immediately the following characterization of the points of 
$(D^2)^\perp\backslash(D^3)^\perp$ in which the germ of 
corresponding reduced Jacobi curve has a constant weight:

\begin{prop}
\label{constJprop} For any $\lambda\in 
(D^2)^\perp\backslash(D^3)^\perp$ the following relation 
holds 
\begin{equation}
\label{difdimJ} {\rm dim}\, {\cal J}^{(i)}(\lambda)- {\rm 
dim}\, {\cal J}^{(i-1)}(\lambda)\leq 1,\quad \forall\, 
i=1,\ldots n-3 
\end{equation}
The germ of reduced Jacobi curve $\widetilde J_\gamma$ at 
$\lambda\in \gamma$ has constant weight iff 
\begin{equation}
\label{Jconstw} {\rm dim}\, {\cal J}^{(n-3)}(\lambda)=2n-4 
\end{equation}
In this case the weight is equal to $(n-3)^2$. 
\end{prop}

\begin{remark}
{\rm From (\ref{difdimJ}) it follows that (\ref{Jconstw}) 
is equivalent to the following relations 

\begin{equation}
\label{Jconstw1} {\rm dim}\, {\cal J}^{(i)}(\lambda)=n-1+i, 
\quad \forall\, i=1,\ldots n-3 
\end{equation}}
\end{remark}
 Denote by ${\cal R}_D$ the set of all 
$\lambda\in (D^2)^\perp \backslash (D^3)^\perp$ such that 
the germ of reduced Jacobi curve $\widetilde J_\gamma$ at 
$\lambda\in \gamma$ has constant weight. By the previous 
proposition ,
\begin{equation}
\label{Jconstw2} {\cal R}_D=\{\lambda\in (D^2)^\perp 
\backslash (D^3)^\perp: {\rm dim}\, {\cal 
J}^{(n-3)}(\lambda)=2n-4 \}. 
\end{equation}
Also $\forall q\in M$ let 

\begin{equation}
\label{calRDq}
{\cal R}_D(q)={\cal R}_D\cap 
T_q^*M
\end{equation}
 and $(D^2)^\perp(q)$ be as in (\ref{D2q}). 
 The question is whether for generic germ of rank 
 $2$-distribution at $q$ the set ${\cal R}_D(q)$ is not 
 empty so that we can apply the theory, presented in subsection 2.1.

For this first we will investigate the following question: 
suppose that the reduced Jacobi curve of the regular 
abnormal extremal $\gamma$ has constant weight; what can be 
said about the corresponding abnormal trajectory 
${\xi}=\pi(\gamma)$? Take some basis $(X_1,X_2)$ in a 
neighborhood of the curve $\xi$ such that $\xi$ is tangent 
to the line distribution spanned by $X_1$ (since our 
considerations are local we always can do it, restricting 
ourselves, if necessary, to some subinterval of $\xi$). For 
any $q\in \xi$ denote by ${\cal T}_\xi^{(i)}(q)$ the 
following subspace of $T_q M$ as follows: 
\begin{equation}
\label {Ti} {\cal T}_\xi^{(i)}(q)={\rm span}\bigl(X_1(q), 
X_2(q),{\rm ad} X_1 (X_2)(q),\ldots, ({\rm ad} 
X_1)^i(X_2)(q)\bigr) 
\end{equation} 
It is easy to see that the subspaces ${\cal T}_\xi^{(i)}(q)$ do 
not depend on the choice of the local basis $(X_1,X_2)$ 
with the above property, but only on the germ of the distribution 
$D$ and the curve $\xi$ at $q$. The property of the curve 
$\xi$ to be abnormal trajectory can be described in terms 
of ${\cal T}_\xi^{(i)}(q)$:

\begin{prop}
\label{JT} If $\gamma$ is abnormal extremal in 
$(D^2)^\perp\backslash (D^3)^\perp$ and $\xi$ is the 
corresponding abnormal trajectory , $\xi=\pi (\gamma)$, 
then $\forall \lambda\in\gamma$ the following relations 
hold 
\begin{equation} \label{TpiJ}{\cal 
T}_\xi^{(i)}\bigl(\pi(\lambda)\bigr)=\pi_*{\cal 
J}^{(i)}(\lambda),
\end{equation}
\begin{equation}
\label{dimTJ} \dim {\cal T}_\xi^{(i)} 
(\pi(\lambda)\bigr)=\dim\,{\cal J}^{(i)}(\lambda)-(n-3).
\end{equation}
\end{prop} 

{\bf Proof.} Let, as before, $H$ be some vector field 
without stationary points 
 tangent to characteristic $1$-foliation $Ab_D$ in a 
neighborhood of $\gamma$. Also, let $\widetilde {\cal X}$ 
be some vector field in a neighborhood of $\gamma$ such 
that $\pi_*\Bigl({\rm span}\bigl(H(\lambda), \widetilde 
{\cal X}(\lambda)\bigr)\Bigr)=D\bigl(\pi(\lambda)\bigr)$.
Then from construction of ${\cal J}^{(i)}$ and relations 
(\ref{part6})-(\ref{hX}) it follows easily that 
\begin{equation}
\label{JadX}
 {\cal J}^{(i)}(\lambda)={\rm span}\Bigl(T_\lambda 
 \Bigl((D^2)^\perp\bigl(\pi(\lambda)\bigr)\Bigr),H(\lambda),
 \widetilde{\cal X}(\lambda),
 {\rm 
ad} H(\widetilde{\cal X})(q),\ldots, ({\rm ad} 
H)^i(\widetilde{\cal X})(q)\Bigr). \end{equation} Take some 
$n$-dimensional submanifold $\Sigma$ of $(D^2)^\perp$, 
passing through $\gamma$ transversal to the fibers 
$(D^2)^\perp\bigl(\pi(\lambda)\bigr)$ for any 
$\lambda\in\gamma$. By construction, $\pi$ projects some 
neighborhood $\widetilde\Sigma$ of $\gamma$ in $\Sigma$ 
bijectively to some neighborhood $V$ of $\xi$ in $M$. 
Taking $$X_1\bigl(\pi(\lambda)\bigr)=\pi_* H (\lambda),\,\, 
X_2 (\pi(\lambda)\bigr)=\pi_*\widetilde{\cal X} 
(\lambda),\quad \forall \lambda\in\widetilde\Sigma$$ and 
using equations (\ref{Ti}), (\ref{JadX}), one obtains 
(\ref{TpiJ}). Relation (\ref{dimTJ}) follows from 
(\ref{TpiJ}) and the fact that 
 the fiber $(D^2)^\perp(q)$ is $(n-3)$-dimensional. This 
 concludes the proof. $\Box$
\medskip

\begin{cor}
\label{Tconstw} The reduced Jacobi curve of the regular 
abnormal extremal $\gamma$ has constant weight iff 
\begin{equation} \label{Tn-3} 
  \dim {\cal T}_\xi^{(n-3)}(q)=n-1,\quad\forall q\in\xi, 
 \end{equation}
 where $\xi=\pi(\gamma)$ is the abnormal trajectory corresponding to 
 $\gamma$.
 \end{cor}
 \begin{remark}
 \label{remreg}{\rm  Note that a smooth curve $\xi$ in $M$, satisfying 
 (\ref{Tn-3}) together with the following relation
  \begin{equation} 
\label{Tn-2} 
  \dim {\cal T}_\xi^{(n-2)}(q)=n-1,\quad\forall q\in\xi, 
 \end{equation}
 is corank 1 abnormal trajectory (see Remark \ref{PMP} for definition of 
 corank). If in addition to (\ref{Tn-3}) and (\ref{Tn-2}) the following relation holds
 \begin{equation}
 \label{regtraj}
  {\cal T}_\xi^{(n-3)}(q)+D^3(q)=T_qM,\quad\forall q\in\xi,
  \end{equation}
  then the curve $\xi$ is regular abnormal extremal.
  In term of a local basis $(X_1,X_2)$ in the 
such that $\xi$ is tangent to the line distribution spanned 
by $X_1$, the condition (\ref{regtraj}) is equivalent to 
the fact that $\forall q\in\xi$ the following condition 
holds 
\begin{equation} \label{regtraj1} {\rm span}\bigl(X_1(q), 
X_2(q),{\rm ad} X_1 (X_2)(q),\ldots, ({\rm ad} 
X_1)^{n-3}(X_2)(q),[X_2,[X_1,X_2]](q)\bigr)=T_qM. 
\end{equation} 
  The assertions of this remark can be deduced without 
 difficulties from the fact that abnormal trajectories are 
 critical points of certain endpoint mapping ( or 
 time$\times$ input/state mapping) and from the expression 
 for the first differential for this mapping (one can use , 
 for example, \cite{agrsar}, section 4).
  }$\Box$
\end{remark}
  \begin{remark} {\rm If the germ of regular abnormal trajectory
   $\xi$ at some 
 point $q_0$
  has corank 1, then the set of $q\in \xi$, satisfying 
  (\ref{Tn-3}),
is open and dense set in some neighborhood of $q_0$ in $\xi$.   
 } $\Box$ 
  \end{remark}

Now we are ready to prove the following genericity result:
   
\begin{prop}
\label{genconct} For generic germ of $(2,n)$-distribution 
 $D$ at $q_0\in M$ ($n\geq 4$) the set ${\cal R}_D(q_0)$, defined in 
(\ref{calRDq}), is a 
nonempty open set in Zariski topology on the linear space 
$(D^2)^\perp(q_0)$, i.e., ${\cal R}_D(q_0)$ is a complement 
to some proper algebraic variety of $(D^2)^\perp(q_0)$. 
\end{prop}

{\bf Proof.} First note that the set 
$(D^2)^\perp(q_0)\backslash {\cal R}_D(q_0)$ is an 
algebraic variety in the linear space $(D^2)^\perp(q_0)$. 
Indeed, choose again a local frame $\{X_i\}_{i=1}^n$ on $M$ 
such that $X_1,X_2$ constitute a local basis of $D$ and 
$X_3$, $X_4$, $X_5$ satisfy (\ref{x345}). Then from 
(\ref{plambda}), definitions of subspaces ${\cal 
J}^{(i)}(\lambda)$ and vector field $\vec h_{X_1,X_2}$ it 
follows that as a basis of spaces ${\cal J}^{(i)}(\lambda)$ 
one can take some vector fields, which are linear 
combination of the fields $\vec u_k$, $\partial_{u_l}$ with 
polynomial in $u_j$ coefficients (here $k,l=1,\ldots,n$, 
$j=4,\ldots,n$). Therefore the set 
\begin{equation}\label{singset} (D^2)^\perp(q_0)\backslash {\cal 
R}_D(q_0)=\{\lambda\in (D^2)^\perp(q_0): {\cal 
J}^{(n-3)}(\lambda)<2n-4\}
\end{equation}
 can be represented as a zero 
level set of some polynomial in $u_j$, $j =4,\ldots,n$. 

Further the coefficients of this polynomial are some 
polynomials in the space of $l_n$-jets of $(2,n)$- 
distributions for some natural $l_n$. We will denote this 
space by ${\rm Jet}_{2,n}(l_n)$. It implies that there 
exists an open set ${\cal U}_n$ in Zariski topology of 
 ${\rm Jet}_{2,n}(l_n)$ 
such that the set ${\cal R}_D(q_0)$ is not empty iff the 
$l_n$-jet of $D$ at $q_0$ belongs to ${\cal U}_n$. Note 
that if the set ${\cal U}_n$ is not empty, then it is dense 
in ${\rm Jet}_{2,n}(l_n)$. 
 Therefore in order to prove our 
proposition it is sufficient to give an example of germ of 
$(2,n)$-distribution such that ${\cal R}_D$ is nonempty. As 
such example one can take distribution $D_0$ spanned by the 
following vector fields 
\begin{equation} 
\label{examprd1} X_1={\partial\over\partial x_1}, \quad 
X_2= {\partial\over\partial x_2}+\sum_{i=1}^{n-3} 
\frac{x_1^i}{i!}{\partial\over\partial 
x_{i+2}}+x_1x_2{\partial\over\partial x_n}, 
\end{equation} 
where $(x_1,\ldots,x_n)$ are some local coordinates on $M$, 
$q_0$=(0,\ldots,0). Using Remark \ref{remreg} and Corollary 
\ref{Tconstw}, it is easy to see that the curve 
$(x_1,0,\ldots,0)$ is regular abnormal trajectory and its 
lift has the reduced Jacobi curve of constant weight. This 
implies that ${\cal R}_{D_0}(q_0)\neq \emptyset$. $\Box$ 

\medskip

Below we give an explicit description of the set ${\cal 
R}_D$ for $n=4,\, 5$ and $6$. In the case $n=4$, small 
growth vector $(2,3,4)$, from (\ref{hX}) it follows 
immediately that ${\cal 
R}_D=(D^2)^\perp\backslash(D^3)^\perp$. A similar result 
holds in the case $n=5$: 
 
\begin{prop}
\label{5weiconst}
For $(2,5)$-distribution with small growth vector $(2,3,5)$ the following 
relation holds
\begin{equation}
\label{XiD5}
 {\cal R}_D=(D^2)^\perp\backslash(D^3)^\perp.
\end{equation}
\end{prop}
{\bf Proof.} Let vector fields $\vec h_{X_1,X_2}$ and $F$ be as in 
(\ref{ham25}) and (\ref{F25}) respectively. Then, using (\ref{j1}), one 
can obtain  by direct computations that
\begin{equation}
\label{hF} [\vec h_{X_1,X_2},F]=u_4 \vec u_5-u_5 \vec 
u_4\quad ({\rm mod}\,\, J^{(1)}) 
\end{equation}
(actually this formula holds for all $n\geq 5$). Hence 
${\rm dim}\, J^{(2)}(\lambda)={\rm dim}\, 
J^{(1)}(\lambda)+1=6$ for all $\lambda \in 
(D^2)^\perp\backslash(D^3)^\perp$, which implies 
(\ref{XiD5}).$\Box$ 

\begin{remark}
\label{D3=5rem}{\rm Let $D$ be $(2,n)$-distribution ($n\geq 
5$) such that ${\rm dim} D^3(q)=4$ for any $q$ in some 
neighborhood $U$. Then from (\ref{hF}) it follows easily 
that $J^{(2)}(\lambda) =J^{(1)}(\lambda)$ for any $q\in 
\pi(U)$. It implies that ${\cal R}_D(q) =\emptyset$ for any 
such $q$ and the theory of subsection 2.1 cannot be applied 
for the reduced Jacobi curves.}$\Box$ 
\end{remark}

 In the case of $(2,6)$-distribution $D$ with growth vector of the type 
 $(2,3,5,\ldots)$  the set ${\cal R}_D$ can be described as follows: Take some $\bar\lambda=(\bar p, q)\in
(D^3)^\perp(q)$ and some vector $v\in D(q)$. Let ${\nu}$ be 
some vector field tangent to $D$ such that $\nu(q)=v$. 
Also, let $(X_1,X_2)$ be a local basis of distribution. 
Then it is easy to see that the number 
 $\bar p\cdot\Bigl[\nu,\bigl[\nu,[X_1,X_2]\bigr]\Bigr](q)$ does not depends on the choice of the vector field $\nu$, so one has a quadratic form 
\begin{equation}
\label{Qn6} v\mapsto Q_{\bar\lambda, 
X_1,X_2}(v)\stackrel{def}{=} \bar 
p\cdot\Bigl[\nu,\bigl[\nu,[X_1,X_2]\bigr]\Bigr](q) 
\end{equation}
on $D(q)$.
 Besides, a change of the local basis of distribution causes to multiplication of this quadratic form on a nonzero constant (which is equal to the determinant 
of the transition matrix between the bases). For 
$(2,6)$-distribution the linear space $(D^3)^\perp(q)$ is 
one-dimensional. Therefore the zero level set 
\begin{equation}
\label{Kn6} {\cal K}(q)=\{v\in D(q): Q_{\bar\lambda, 
X_1,X_2}(v)=0\} 
\end{equation}
of $Q_{\bar\lambda, X_1,X_2}$ is the same for all 
$\bar\lambda\in (D^3)^\perp(q)\backslash(0,q)$ and any 
local basis $X_1,X_2$ of the distribution. 
\begin{prop}
\label{RD6prop} For $(2,6)$-distribution $D$ with the small 
growth vector of the type $(2,3,5,\ldots)$ the following 
relation holds 
\begin{equation}
\label{XiD6} {\cal R}_D(q)= \{\lambda\in 
(D^2)^\perp(q):\pi_*\bigl({\rm ker}\,  
\sigma|_{(D^2)^\perp}(\lambda)\bigr)\notin {\cal K}(q)\}. 
\end{equation}
The set ${\cal R}_D(q)\neq \emptyset$ iff the small growth 
vector of $D$ at $q$ is equal to $(2,3,5,6)$. 
\end{prop}

{\bf Proof.} As before, complete some basis $X_1,X_2$ of $D$ 
to the frame $\{X_i\}_{i=1}^6$ on $M$ such that $X_3$, 
$X_4$, $X_5$ satisfy (\ref{x345}). Let $c_{ji}^k$ be the 
structural functions of this frame, i.e., the functions, 
satisfying $[X_i,X_j]=\sum_{k=1}^6 c_{ji}^k X_k$. Then from 
(\ref {hF}) by straightforward calculation it follows that 
\begin{equation}
\label{HJ3} [\vec h_{X_1,X_2},u_4 \vec u_5-u_5 \vec u_4]= 
[u_4 \vec u_2-u_5 \vec u_1,u_4 \vec u_5-u_5 \vec u_4]= 
\alpha_6 \vec u_6 \,\,\bigl({\rm mod}\,\,{\rm span}(\vec 
u_4,\vec u_5, J^{(1)})\bigr), \end{equation} where 
\begin{equation} \label{a6mal} 
\alpha_6=c_{52}^6u_4^2-(c_{42}^6+c_{51}^6)u_4u_5+ c_{41}^6 
u_5^2 
\end{equation}
From (\ref{hF}) and (\ref{HJ3}) it follows that if 
$J^{(3)}(\lambda)=J^{(2)}(\lambda)$, then $\alpha_6=0$. 
Conversely, if $\alpha_6=0$ then by (\ref{HJ3}) 
$J^{(3)}(\lambda)\subset {\rm span}(\vec u_4,\vec u_5, 
J^{(1)})$. But by construction $ J^{(3)}(\lambda)\subset 
\vec e(\lambda)^\angle$ (where $\vec e(\lambda)^\angle$ is 
the skew-symmetric complement of $\vec e(\lambda)$ in 
$T_\lambda T^*M$). This together with (\ref{hF}) implies 
that $J^{(3)}(\lambda)=J^{(2)}(\lambda)$. So, $\lambda\in 
{\cal R}_D(q)$ iff $\alpha_6(\lambda)\neq 0$. To prove 
(\ref{XiD6}) it remains to note that $$Q_{\bar\lambda, 
X_1,X_2}(\pi_*\bigl(\vec 
h_{X_1,X_2}(\lambda)\bigr)=C\alpha_6(\lambda),$$ where $C$ 
is a nonzero constant. 
 The last 
assertion of the proposition follows from the fact that 
$\alpha_6\equiv 0$ iff $c_{ji}^6=0$,where $i=1,2$, $j=4,5$, 
or, equivalently, that $\dim\, D^4(q)=5$. $\Box$ 
\medskip

{\bf 2.3 Fundamental form of distribution and its 
properties.} For any $\lambda\in {\cal R}_D$ take 
characteristic curve $\gamma$, passing through $\lambda$. 
Let ${\cal A}_\lambda$ be the fundamental form of the 
reduced Jacobi curve $\widetilde J_\gamma$ of $\gamma$ at 
$\lambda$. By construction ${\cal A}_\lambda$ is degree $4$ 
homogeneous function on the tangent line to $\gamma$ at 
$\lambda$. In the previous subsection to any (local) basis 
$(X_1, X_2)$ of distribution $D$ we assigned the vector 
field $\vec h_{_{X_1,X_2}}$ tangent to characteristic 
$1$-foliation $Ab_D$ (see (\ref{ham25})). Let 
\begin{equation}
\label{densX12} A_{_{X_1,X_2}}(\lambda)= {\cal 
A}_\lambda(\vec h_{_{X_1,X_2}}(\lambda)) \end{equation} In 
this way 
 to any (local) basis 
$(X_1, X_2)$ of distribution $D$ we assign the function $ 
A_{_{X_1,X_2}}$ on ${\cal R}_D$. 

\begin{remark}
\label{densrem} {\rm If we consider parametrization 
$t\mapsto \widetilde J_\gamma(e^{t\vec 
h_{_{X_1,X_2}}}\lambda)$ of the reduced Jacobi curve of 
$\gamma$, then $A_{_{X_1,X_2}}(\lambda)$ is the density of 
fundamental form of this curve w.r.t. parametrization $t$ 
at $t=0$.} $\Box$ 
\end{remark}

Let $\tilde X_1,\tilde X_2$ be another basis of the 
distribution $D$. Then there exist functions 
$\{\nu_{ij}\}_{i,j=1}^2$ on $M$ such that 

$$
\begin{array}{l}
\tilde {X_1}=\nu_{11} X_1+\nu_{12} \,X_2 \\ \tilde 
{X_2}=\nu_{21}\,X_1+\nu_{22}\, X_2 
\end{array}$$
By direct computation one has
\begin{equation}
\label{tildehh} \vec h_{_{\tilde X_1,\tilde 
X_2}}(\lambda)=\Delta^2\bigl(\pi(\lambda)\bigr) \vec 
h_{_{X_1,X_2}}(\lambda), 
\end{equation}
where $\Delta$ is equal to determinant of transition matrix 
from the basis $(X_1,X_2)$ to the basis $(\tilde X_1,\tilde 
X_2)$, i.e., $\Delta=\nu_{11}\nu_{22}-\nu_{12}\nu_{21}$. 
From this and homogeneity of ${\cal A}$ it follows that 

\begin{equation}
\label{transAdist} 
 A_{_{\tilde X_1,\tilde X_2}}(\lambda)=
\Delta\bigr(\pi(\lambda)\bigl)^8 A_{_ {X_1,X_2}}(\lambda) 
\end{equation}

Therefore for any $q\in M$ such that ${\cal R}_D(q)\neq 
\emptyset$ the restriction of $A_{_ {X_1,X_2}}$ to ${\cal 
R}_D(q)$ is well defined function, up to multiplication on 
positive constant, or well defined element of "positive 
projectivization" of the space of the functions on ${\cal 
R}_D(q)$. We will call it {\it fundamental form of the rank 
2 distribution D at the point $q$}. From now on we will 
write $\vec h$ instead of $\vec h_{_{X_1,X_2}}$ and $A$ 
instead of $A_{_{X_1,X_2}}$ without special mentioning. 

\begin{remark}
\label{bryantps}
 {\rm According to subsection 2.1 ( see the 
sentence after formula (\ref{sch})) any abnormal extremals 
of $(2,n)$-distribution $D$ lying in ${\cal R}_D$ carries 
the canonical projective structure. It can be shown that in 
the case $n=4$, small growth vector $(2,3,4)$, our 
canonical projective structure defined on abnormal 
extremals ( and therefore also on abnormal trajectories) 
coincides with the projective structure on abnormal 
trajectories, introduced in \cite{bryant} (see Proposition 
5 there).
Note also that by Remark \ref{m=1A=0} and relation 
(\ref{dimjac}) in the case $n=4$ the fundamental form is 
identically equal to zero.} 
\end{remark}
\begin{remark}
\label{remsmooth} {\rm Using Remark \ref{ratrem} it is easy 
to see that the fundamental form $A(\lambda)$ is a smooth 
function for all $\lambda\in {\cal R}_D$: one can choose 
the coordinate representation of the curves $t\mapsto 
\widetilde J_\gamma(e^{t\vec h}\lambda)$ smoothly depending 
on $\lambda$ and use the fact that the operation of 
differentiation by $t$ in coordinates corresponds to the 
operation ${\rm ad}\, \vec h$ due to the relation 
(\ref{addt}).} 
\end{remark}

In fact one can say much more about algebraic structure of 
the fundamental form of distribution. 
\begin{prop}
\label{rathom4} For any $q\in M$ such that ${\cal 
R}_D(q)\neq \emptyset$ the fundamental form of the rank $2$ 
distribution $D$ at the point $q$ is degree $4$ homogeneous 
rational function on $(D^2)^\perp(q)$, defined up to 
multiplication on positive constant. \end{prop} 

{\bf Proof.} First let us prove that the fundamental form 
at $q$ is rational function on $(D^2)^\perp(q)$. From 
Remark \ref{ratrem} it follows that in order to do this it 
is sufficient to show that the parametrized reduced Jacobi 
curves $t\mapsto \widetilde J_\gamma(e^{t\vec h}\lambda)$ 
have coordinate representations $t\mapsto S_\lambda(t)$ 
such that for any natural $l$ all entries of 
$S_\lambda^{(l)}(0)$, as functions of $\lambda$, are 
rational functions on the fibers $(D^2)^\perp(q)$. For this 
choose the following 
 $(2n-3)$ vector fields   on $(D^2)^\perp$:
\begin{equation}
\label{basd2perp}
\partial_\theta, X,\partial_{u_6},\ldots,\partial_{u_n},F, Y_4,\ldots Y_{n-1}, Z,
\vec e,\vec h,
\end{equation}
where
\begin{eqnarray}
~&~ Y_k =\displaystyle{u_{k+1}\vec {u_k}-u_k\vec
{u_{k+1}}+\sum_{i=1}^3 \Bigl(u_{k+1}\{u_i,u_k\}-
u_k\{u_i,u_{k+1}\}\Bigr)\partial_{u_i}},\label{vecYk}\\
~&~Z=\displaystyle{u_4\vec {u_5}+u_5 \vec {u_4}+\sum_{i=1}^3
\Bigl(u_4\{u_i,u_5\}+u_4\{u_i,u_4\}\Bigr)\partial_{u_i}}\label{vecz}
\end{eqnarray}
(here $\{u_i,u_j\}$ are Poisson brackets of the Hamiltonians $u_i$
and $u_j$, i.e.,  $\{u_i,u_j\}=du_j(\vec{u_i})$). 
Let 
\begin{equation}
\label{wl}
W_\lambda=\Bigl(\vec e(\lambda)^\angle\cap T_\lambda 
(D^2)^\perp\Bigr)/{\rm span} \bigr({\rm ker}\, \sigma|_{(D^2)^\perp}(\lambda) ,
\vec e(\lambda)\bigl)
\end{equation}
(here by $\vec e(\lambda)^\angle$ we mean a skew-symmetric 
complement of $\vec e(\lambda)$ in $T_\lambda T^*M$). Then 
under identification (\ref{tanid}) the reduced Jacobi curve 
$\widetilde J_\gamma$ lives in Lagrange Grassmannina 
$L(W_\lambda)$ of symplectic space $W_\lambda$. Denote by 
${\cal P}$ the set of all $\lambda\in (D^2)^\perp$ such 
that the vector fields (\ref{basd2perp}) at $\lambda$ 
constitute a basis of $T_\lambda(D^2)^\perp$. Evidently, 
for any $q\in M$ the set ${\cal P}\cap (D^2)^\perp(q)$ is a 
nonempty open set in Zariski topology on the linear space 
$(D^2)^\perp(q)$. For any $\lambda\in (D^2)^\perp$ the 
first $2(n-3)$ vectors in (\ref{basd2perp}) belong to $\vec 
e(\lambda)^\angle$. Therefore, for any $\lambda\in {\cal 
P}$ the images of the first $2(n-3)$ vectors in 
(\ref{basd2perp}) under the canonical projection from 
$\Bigl(\vec e(\lambda)^\angle\cap T_\lambda 
(D^2)^\perp\Bigr)$ to $W_\lambda$ constitute the basis of 
the space $W_\lambda$. Introduce in $W_\lambda$ the 
coordinates w.r.t. this basis and suppose that $t\mapsto 
S_\lambda(t)$ is the corresponding coordinate 
representation of the curve $t\mapsto \widetilde J_\gamma(e^{t\vec h}\lambda)$,
 $\widetilde J_\gamma(e^{t\vec h}\lambda)=\{x, S_\lambda(t)):x\in 
\mathbb{R}^{n-3}\}$. Then from (\ref{ident}) and 
(\ref{addt}) it follows that for any natural $l$ all 
entries of the matrix $S^{(l)}_\lambda(0)$ are some 
rational combinations of some coordinates of the vectors of 
the type 
\begin{equation}
\label{adn}
({\rm ad}\vec h)^j(\partial_\theta)(\lambda),\,\,\,
({\rm ad}\vec h)^j( X)(\lambda),\,\,\,
({\rm ad}\vec h)^j(\partial_{u_i})(\lambda),\,\,\,\,\,6\leq i\leq n,
\,\,\,\,\, 1 \leq j \leq l.
\end{equation}
 w.r.t. the basis (\ref{basd2perp}).
But from the form  of the vector fields $Y_i$  and $Z$ it is clear
that coordinates
 of the vectors from (\ref{adn})  w.r.t. the basis (\ref{basd2perp})
 are rational functions on the fibers
$(D^2)^\perp(q)$. So, for any $q$ the fundamental form at 
$q$ is a rational function on the fiber $(D^2)^\perp(q)$. 

 Now let us show that the fundamental form is homogeneous of degree $4$.
Indeed, it is clear that $$\delta_{a\,*}{\cal 
J}(\lambda)=J(\delta_a(\lambda)),$$ where $\delta_a$ is 
the homothety defined by (\ref{homoth}). This together with 
Remark \ref{homrem} implies that $\delta_{a\,*}$ induces 
the symplectic transformation from $W_\lambda$ to 
$W_{\delta _a(\lambda)}$, which transforms 
the curve $\widetilde J_\gamma$ to the curve $\widetilde 
J_{\delta_a(\gamma)}$. Therefore the following identity 
holds 
\begin{equation}
\label{covar} {\cal A}_{\delta_a(\lambda)}\bigr(\delta
_{a\,*}\vec h(\lambda)\bigl)= {\cal A}_\lambda\bigr(\vec 
h(\lambda)\bigl). 
\end{equation}
On the other hand, one has $$\vec h \bigr(\delta
_a(\lambda)\bigl)=a\,\delta_{a\,*}\vec h(\lambda).$$ Hence 
$$A\bigl(\delta_a \lambda)={\cal A}_{\delta 
_a(\lambda)}\Bigl(\vec h \bigr(\delta 
_a(\lambda)\bigl)\Bigr)=a^4 {\cal A}_{\delta 
_a(\lambda)}\bigl(\delta_{a\,*}\vec 
h(\lambda)\bigr)=a^4\,{\cal A}_\lambda\bigr(\vec 
h(\lambda)\bigl)=a^4\,A(\lambda).$$ So $A$ is homogeneous 
of degree $4$. $\Box$

In the case $n=5$ and small growth vector $(2,3,5)$ one can 
look at the fundamental form of the distribution $D$ from 
the different point of view. In this case (in contrast to 
generic $(2,n)$-distributions with $n>5$) there is only one 
abnormal trajectory starting at given point $q\in M$ in 
given direction (tangent to $D(q)$). All lifts of this 
abnormal trajectory can be obtained one from another by 
homothety. So they have the same, up to symplectic 
transformation, Jacobi curve. It means that one can 
consider Jacobi curve and fundamental form of this curve on 
abnormal trajectory instead of abnormal extremal. 
Therefore, to any $q\in M$ one can assign a homogeneous 
degree $4$ rational function $\AA_q$ on the plane $D(q)$ in 
the following way: 
\begin{equation}
\label{tangfund} \AA_q(v)\stackrel{def}{=}{\cal 
A}_\lambda(H) 
\end{equation} for any $v\in D(q)$, where 
\begin{equation}
\label{condlh} \pi(\lambda)=q,\quad\pi_*H=v,\quad H\in {\rm 
ker}\,\sigma|_{(D^2)^\perp}(\lambda). \end{equation} and 
the righthand side of (\ref{tangfund}) is the same for any 
choice of $\lambda$ and $H$, satisfying (\ref{condlh}). 
$\AA_q$ will be called {\it tangential fundamental form} of 
the distribution $D$ at the point $q$. We stress that the 
tangential fundamental form is the well defined function on 
$D(q)$ and not the class of functions under positive 
projectivization. 

The analysis of the algebraic structure presented in the 
proof of Proposition \ref{rathom4} is rather rough. In the 
sequel we will show that for $n=5$ fundamental form is 
always polynomial on $(D^2)^\perp(q)$ (defined up to 
multiplication on a positive constant), while for $n>5$ it 
is nonpolynomial rational function for generic distribution. 

{\bf 2.4 Projective curvature of rank 2 distribution with 
nonzero fundamental form.}  
Denote by 
\begin{equation}
\label{aneq0} \aleph_D=\{\lambda\in 
{\cal R}_D: {\cal A}_\lambda\neq 0\} 
\end{equation}
Suppose that the set $\aleph_D$ is not empty. 

\begin{remark}{\rm 
For $n=5$ the set $\aleph_D$ is empty iff distribution is 
locally equivalent to so-called free nilpotent 
$(2,5)$-distribution (see Example 1 and Remark 
\ref{nilpcart} in section 3). Our conjecture is that 
$\aleph_D$ is empty iff $n$ is equal to the dimension of 
the free nilpotent $r$-step Lie algebra $g_{r,2}$ with two 
generators 
 for some $r\geq 3$ and $D$ is locally equivalent 
to the left-invariant distribution on a Lie group with Lie 
algebra $g_{r,2}$ such that this distribution is spanned by 
the generators of $g_{r,2}$}.$\Box$ 
\end{remark}

As was mentioned in subsection 2.1, for any characteristic 
leaf lying in $\aleph_D$, the corresponding Jacobi curve 
has canonical parameter. In other words, any such leaf has 
canonical parameter. 
Define the vector field $\overrightarrow {h_A}$ on 
$\aleph_D$ such that its integral curves are the 
characteristic leaves parameterized by their canonical 
parameter (i.e., the field $\overrightarrow {h_A}$ is given 
by velocities of the characteristic leaves parameterized by 
their canonical parameter). By construction, the vector 
field $\overrightarrow {h_A}$ is invariant of the 
distribution $D$. We call it {\it the canonical Hamiltonian 
vector field } of distribution $D$. 
Further, for any $\lambda\in\aleph_D$ take the 
characteristic leaf $\gamma$, passing through $\lambda$. 
Denote by $\rho{_{_D}}(\lambda)$ the projective Ricci 
curvature of the reduced Jacobi curve $\tilde J_\gamma$ 
at $\lambda$. In other words, $\rho{_{_D}}(\lambda)$ is 
equal to the Ricci curvature of the curve $\tau\mapsto 
\widetilde J_\gamma(e^{\tau \overrightarrow{h_A}}\lambda)$ 
at the point $\tau=0$. 
So, to given rank 2 distribution $D$ we assign canonically the 
function 
\begin{equation}
\label{funinv2n} \rho_{_{D}}: \aleph_D\mapsto\mathbb{R} 
\end{equation}
This function is called {\it a projective Ricci curvature} 
of distribution $D$. 

Now we give a method for computation of projective 
curvature $\rho{_{_D}}$. Take some local basis $X_1,X_2$ of 
$D$. Let again $\vec h=\vec h_{X_1 X_2}$ and 
$A=A_{X_1,X_2}$ be as in 
(\ref{ham25}), and (\ref{densX12}) 
respectively. 
Also denote by $\rho(\lambda)$ the Ricci curvature of the 
parameterized curve $t\mapsto\widetilde J_\gamma(e^{t \vec 
h}\lambda)$ at the point $t=0$. Note that in contrast to 
$\rho{_{_D}}(\lambda)$, the function $\rho(\lambda)$ 
certainly depends on the local basis of distribution. Using 
the reparametrization rule (\ref{rhorep}) for Ricci 
curvature, one can easily express the projective curvature 
$\rho{_{_D}}(\lambda)$ by $\rho(\lambda)$ and $A(\lambda)$. 
Indeed, let $\tau$ be the canonical parameter on $\gamma$ 
and $t$ be parameter defined by the field $\vec h$. Then by 
Remark \ref{denscan} 
\begin{equation}
\label{normarb} d\tau= 
\sqrt[4]{|A(e^{t\vec h}\lambda)|}\,\, dt 
\end{equation}
Suppose that $t=\vf(\tau)$. Then by (\ref{normarb}) 
\begin{equation}
\label{vf'normarb} \vf'(\tau)=\frac{1}{ 
\sqrt[4]{|A(e^{t\vec h}\lambda)|}} 
\end{equation}

Recall that the Jacobi curves under consideration have the 
weight equal to $(n-3)^2$. So, by (\ref{rhorep})
\begin{equation}
\label{transrhodist} 
\rho{_{_D}}\left(e^{\tau\vec{h_A}}\lambda\right)= 
\rho\left(e^{\vf(\tau)\vec h}\lambda\right)\bigl( 
\vf'(\tau)\bigr)^2+\frac{(n-3)^2}{3}\mathbb{S}\bigl(\vf(\tau)\bigr), 
\end{equation}
where $\mathbb {S}(\vf)$ is Schwarzian of the function 
$\vf$, defined by (\ref{sch}). One can check that 
Schwarzian satisfies the following relation 
\begin{equation}
\label{schwarzsimple} \mathbb 
{S}\bigl(\vf(\tau)\bigr)=-\frac{y''(\tau)}{y(\tau)}, 
\end{equation}
where $$y(\tau)=\frac{1}{\sqrt{\vf'(\tau)}}.$$ By (\ref{vf'normarb}), 
$$y(\tau)=\sqrt[8]{|A(e^{\tau\overrightarrow {h_A}}\lambda)|}$$ 
Substituting this in (\ref{schwarzsimple}) and using (\ref{normarb}) 
we obtain 

$$\mathbb {S}\bigl(\vf(\tau)\bigr)= -\frac{\frac{d^2}{d\tau^2} 
\Bigl(\sqrt[8]{|A(e^{\tau\overrightarrow {h_A}}\lambda)|}\Bigr)} 
{\sqrt[8]{|A(e^{\tau\overrightarrow {h_A}}\lambda)|}}=-\frac{1} 
{\sqrt[8]{|A(e^{t\overrightarrow h}\lambda)|^3}}\frac{d}{dt} 
\left(\frac{1}{\sqrt[4]{|A(e^{t\vec h}\lambda)|}} 
\frac{d}{dt}\left(\sqrt[8]{|A(e^{t\vec h}\lambda)|}\right)\right)=$$ 
\begin{equation}
\label{yvfA} =\frac{1} {\sqrt[8]{|(A(e^{t\vec h}\lambda)|^3}} 
\frac{d^2}{dt^2}\left(|A(e^{t\vec h}\lambda)|^{-\frac{1}{8}}\right)= 
\frac{\vec h\circ\vec h\left(|A(e^{t\vec 
h}\lambda)|^{-\frac{1}{8}}\right)} {\sqrt[8]{|(A(e^{t\vec 
h}\lambda)|^3}} 
\end{equation}
Finally, substituting (\ref{yvfA}) with $t=0$ in 
(\ref{transrhodist}) we get 

\begin{equation}
\label{projcurv} 
\rho{_{_D}}=\frac{\rho}{\sqrt{|A|}}+\frac{(n-3)^2}{3} 
\frac{\vec h \circ\vec h 
\bigl(|A|^{-\frac{1}{8}})}{\sqrt[8]{|A|^3}}. 
\end{equation}
The last formula can be rewritten also as follows
\begin{equation}
\label{projcurv1} \rho{_{_D}}=\frac{\rho 
A^2-\frac{(n-3)^2}{24}\vec h \circ\vec h (A) A+\frac{3 
(n-3)^2}{64}\bigl(\vec h(A)\bigr)^2}{|A|^{\frac {5}{2}}}. 
\end{equation}
Since $\rho_D$ is well defined function on $\aleph_D$ and 
$A$ is degree $4$ homogeneous rational function on 
$(D^2)^\perp(q)$, defined up to multiplication on a 
positive constant, the numerator
\begin{equation}
\label{secfund} {\cal C}\stackrel{def}{=} \rho 
A^2-\frac{(n-3)^2}{24}\vec h \circ\vec h (A) A+\frac{3 
(n-3)^2}{64}\bigl(\vec h(A)\bigr)^2 \end{equation} of 
(\ref{projcurv1}) is degree $10$ homogeneous function on 
$(D^2)^\perp(q)$, defined up to multiplication on a 
positive constant. This function will be called {\it a 
second fundamental form} of distribution $D$. The second 
fundamental form is rational function on $(D^2)^\perp(q)$, 
because $A$ is rational and also $\rho$ is rational, which 
follows from the same arguments as in Proposition 
\ref{rathom4}. In the case $n=5$ the second fundamental 
form is polynomial, which will follow from Theorem \ref{algstruct} 
below.

\section{Calculation of invariants of (2,5)-distributions } 
\indent \setcounter{equation}{0}
 
In the present section we give explicit formulas for 
computation of the fundamental form and projective Ricci 
curvature in the case of rank 2 distribution on 
$5$-dimensional manifold (as before we assume that the 
small growth vector is $(2,3,5)$). We demonstrate these 
formulas on several examples, showing simultaneously the 
efficiency of our invariants in proving that the rank 2 
distributions are not equivalent. 

{\bf 3.1 Preliminaries.} In order to obtain these formulas 
we need more facts from the theory of curves in 
Grassmannian $G_m(W)$ of half-dimensional subspaces (here 
${\rm dim}\, W=2m$) and in Lagrange Grassmannian $L(W)$ 
w.r.t. to some symplectic form on $W$, developed in 
{\cite{jac1}, \cite{jac2}). Below we present all necessary 
facts from the mentioned papers together with several new 
useful arguments. 

Fix some $\Lambda\in G_m(W)$. As before , let 
$\Lambda^\pitchfork$ be the set of all $m$-dimensional 
subspaces of $W$ transversal to $\Lambda$. Note that any 
$\Delta\in \Lambda^\pitchfork$ can be canonically 
identified with $W/\Lambda$. Keeping in mind this 
identification and taking another subspace $\Gamma\in 
\Lambda^\pitchfork$ one can define the operation of 
subtraction $\Gamma-\Delta$ as follows 
$$\Gamma-\Delta\stackrel{def}{=}\langle\Delta,\Gamma,\Lambda\rangle\in 
{\rm Hom}\,(W/\Lambda,\Lambda).$$ It is clear that the set 
$\Lambda^\pitchfork$ provided with this operation can be 
considered as the affine space over the linear space ${\rm 
Hom}\, (W/\Lambda,\Lambda)$.

Consider now some ample curve $\Lambda(\cdot)$ in $G_m(W)$. 
Fix some parameter $\tau$. By assumptions 
$\Lambda(t)\in\Lambda(\tau)^\pitchfork$ for all $t$ from a 
punctured neighborhood of $\tau$. We obtain the curve 
$t\mapsto\Lambda(t)\in\Lambda(\tau)^\pitchfork$ in the 
affine space $\Lambda(\tau)^\pitchfork$ with the pole at 
$\tau$. We denote by $\Lambda_\tau(t)$ the identical 
embedding of $\Lambda(t)$ in the affine space 
$\Lambda(\tau)^\pitchfork$. First note that the velocity 
${\partial\over\partial t}\Lambda_\tau(t)$ is well defined 
element of ${\rm Hom} (W/\Lambda, \Lambda)$. 
Fixing an ``origin'' in $\Lambda(\tau)^\pitchfork$ we make 
$\Lambda_\tau(t)$ a vector function with values in ${\rm 
Hom}\,(W/\Lambda,\Lambda)$ and with the pole at $t=\tau$. 
Obviously, only free term in the expansion of this function 
to the Laurent series at $\tau$ depends on the choice of 
the ``origin'' we did to identify the affine space with the 
linear one. More precisely, the addition of a vector to the 
``origin'' results in the addition of the same vector to 
the free term in the Laurent expansion. In other words, for 
the Laurent expansion of a curve in an affine space, the 
free term of the expansion is an element of this affine 
space. Denote this element by $\Lambda^0(\tau)$. The curve 
$\tau\mapsto\Lambda^0(\tau)$ is called the {\it derivative 
curve} of $\Lambda(\cdot)$. 
 
If we restrict ourselves to the Lagrange Grassmannian 
$L(W)$, i.e. if all subspaces under consideration are 
Lagrangian w.r.t. some symplectic form $\bar\sigma$ on $W$, 
then from Remark \ref{lagr} it follows easily that the set 
$\Lambda^\pitchfork_L$ of all Lagrange subspaces 
transversal to $\Lambda$ can be considered as the affine 
space over the linear space of all self-adjoint mappings 
from $\Lambda^*$ to $\Lambda$, the velocity 
${\partial\over\partial t}\Lambda_\tau(t)$ is well defined 
self-adjoint mappings from $\Lambda^*$ to $\Lambda$, and 
the derivative curve $\Lambda^0(\cdot)$ consist of Lagrange 
subspaces. Besides if the curve $\Lambda(\cdot)$ is 
nondecreasing rank $1$ curve in $L(W)$, then 
${\partial\over\partial t}\Lambda_\tau(t)$ is a nonpositve 
definite rank $1$ self-adjoint map from $\Lambda^*$ to 
$\Lambda$ and for $t\neq \tau$ there exists a unique, up to 
the sign, vector $w(t,\tau) \in \Lambda(\tau)$ 
  such that for any $v\in \Lambda(\tau)^*$
\begin{equation}
\label{defw1} \langle v,{\partial\over\partial 
t}\Lambda_\tau(t) v\rangle= - \langle v,w(t,\tau)\rangle^2.
\end {equation} The properties of vector function $t\mapsto 
w(t,\tau)$ for a rank $1$ curve of constant weight in 
$L(W)$ can be summarized as follows ( see \cite{jac1}, 
section 7, Proposition 4 and Corollary 2): 

\begin{prop}
\label{constcor} If $\Lambda(\cdot)$ is a rank 1 curve of 
constant weight in L(W), then for any $\tau\in I$ the 
function $t\mapsto w(t,\tau)$ has a pole of order $m$ at 
$t=\tau$. Moreover, if we write down the expansion of 
$t\mapsto w(t,\tau)$ in Laurent series at $t=\tau$, 
$$w(t,\tau)=\sum_{i=1}^{m}e_i(\tau)(t-\tau)^{i-1-l}+O(1),$$ 
then the vector coefficients $e_1(\tau),\ldots,e_m(\tau)$ 
constitute a basis of the subspace $\Lambda(t)$. 
\end{prop}

The basis of the vectors $e_1(\tau),\ldots, e_m(\tau)$ from 
the previous proposition is called a {\it canonical basis} 
of $\Lambda(\tau)$. Further for given $\tau$ take the 
derivative subspace $\Lambda^0(\tau)$ and let 
$f_1(\tau),\ldots, f_m(\tau)$ be a basis of 
$\Lambda^0(\tau)$ dual to the canonical basis of 
$\Lambda(\tau)$, i.e. $\bar\sigma (f_i(\tau), 
e_j(\tau))=\delta_{i,j}$. The basis 
$$(e_1(\tau),\ldots,e_m(\tau),f_1(\tau),\ldots,f_m(\tau))$$ 
of whole symplectic space $W$ is called {\it the canonical 
moving frame} of the curve $\Lambda(\cdot)$. Calculation of 
structural equation for the canonical moving frame is 
another way to obtain symplectic invariants of the curve 
$\Lambda(\cdot)$. 

For the reduced Jacobi curves of abnormal extremals of 
$(2,5)$-distribution $m$ is equal to $2$. So we restrict 
ourselves to this case. For $m=2$ the structural equation 
for the canonical moving frame has the following form (for 
the proof see \cite {jac2} Section 2, Proposition 7): 

\begin{equation}
\label{structeq} \left( 
\begin{array}{l} e_1'(\tau)\\  e_2'(\tau)\\  f_1'(\tau)\\ 
f_2'(\tau)\end{array}\right) = \left(\begin{array}{cccc} 
0&3 &0&0\\ \frac{1}{4} \rho(\tau)&0&0 &4\\ -\left(\frac 
{35}{36} A(\tau) -\frac{1}{8} 
\rho(\tau)^2+\frac{1}{16}\rho''(\tau)\right)& 
-\frac{7}{16}\rho'(\tau) &0&-\frac {1}{4}\rho(\tau) \\ 
-\frac{7}{16} \rho'(\tau) &-\frac{9}{4} \rho(\tau) & - 3 
&0\end{array}\right)\left(\begin{array}{l} 
e_1(\tau)\\e_2(\tau)\\f_1(\tau)\\f_2(\tau)\end{array}\right), 
\end{equation}
where $\rho(\tau)$ and $A(\tau)$ are the Ricci curvature 
and the density of fundamental form of the parametrized 
curve $\Lambda(\tau)$ respectively. 

One can express $e_2(\tau)$ by $e_1'(\tau)$ using the first 
equation of (\ref{structeq}), then $f_2(\tau)$ by 
$e_1(\tau)$ and $e_1''(t)$ using the second equation of 
(\ref{structeq}), then $f_1(\tau)$ by $e_1(\tau)$, 
$e_1'(\tau)$ and $e_1^{(3)}$ using the forth equation of 
(\ref{structeq}). Finally substituting all this to the 
third equation of (\ref{structeq}) one obtains the 
following useful Proposition 

\begin{prop}
Suppose that $\Lambda(t)$ is rank 1 curve of the constant 
weight in $L(W)$ and $e_1(t)$ is the first vector in the 
canonical basis of $\Lambda(t)$ . Then $e_1(t)$ satisfies 
the following relation: 
\begin{equation}
\label{eqe1} e_1^{(4)}=\left(35\,A - \frac{81}{16} 
\rho^2-\frac{9}{4}\rho''\right)e_1 - \frac{15}{2}\rho' 
e_1'-\frac{15}{2}\rho e_1'' 
\end{equation}
\end{prop}

The previous proposition says that in order to find $\rho$ 
and $A$ (which actually constitute a complete system of 
symplectic invariants of the parametrized curve 
$\Lambda(\cdot)$) it is sufficient to know the first vector 
$e_1(\tau)$ in the canonical basis of $\Lambda(\cdot)$. The 
following proposition gives a simple way to find the  
vector $e_1(\tau)$. 
\begin{prop}
\label{canlem1}
 Let $\Lambda(\tau)$ be a rank $1$ nondecreasing curve of constant weight 
 in the Lagrange Grassmannian $L(W)$, where ${\rm dim}\, W=4$. Then the 
 first vector $e_1(\tau)$ of the canonical basis of $\Lambda(\tau)$ 
can be uniquely (up to the sign)  determined by the following two conditions 
 
 \begin{equation}
 \label{dire1}
 \mathbb{R} e_1(\tau)={\cal D}^{(1)}\Lambda(\tau)^\angle,
 \end{equation}
\begin{equation}
\label{seprime}
\bar\sigma(e_1''(\tau),e_1'(\tau))=36,
\end{equation}
 where the subspace  ${\cal D}^{(1)}\Lambda(\tau)$ is as in 
 (\ref{primext}) and ${\cal D}^{(1)}\Lambda(\tau)^\angle$ 
 is its skew-symmetric complement.
 \end{prop}
 
 {\bf Proof.} 
 The relation (\ref{seprime}) follows directly from the first two equations of 
(\ref{structeq}). To prove (\ref{dire1}) note that 
 from (\ref{structeq}) it is clear that
\begin{equation}
\label{d1can}
 {\cal D}^{(1)}\Lambda(\tau)={\rm span}\,
 \bigl(e_1(\tau),e_2(\tau),f_2(\tau)\bigr).
\end{equation}
 But from definition of canonical moving frame it 
follows that $$\Bigl({\rm 
span}\,\bigl(e_1(\tau),e_2(\tau),f_2(\tau)\bigr)\Bigr)^\angle=\mathbb 
{R} e_1(\tau),$$ which together with (\ref{d1can}) implies 
(\ref{dire1}). Finally, the vector $e_1(t)$ is determined 
by (\ref{dire1}) and (\ref{seprime}) uniquely, up to the 
sign: the first relation gives the direction of $e_1(t)$ 
and the second ''normalizes'' this direction. $\Box$ 

\vskip .2in 

{\bf 3.2 Application to $(2,5)$-distributions.} Choose 
some local basis $(X_1,X_2)$ of $(2,5)$-distribution and 
complete it by the fields $X_3$ ,$X_4$, and $X_5$, 
satisfying (\ref{x345}), to the local frame on $M$. Such 
frame $(X_1,X_2,X_3,X_4,X_5)$ will be called {\it adapted} 
to the distribution $D$. If instead of (\ref{x345}) one has 
\begin{equation}
\label{x345s} X_3=[X_1,X_2],\,\,\, 
X_4=[X_1,[X_1,X_2]]=[X_1,X_3],\,\,\,\, 
X_5=[X_2,[X_2,X_1]]=[X_3,X_2], 
\end{equation}
the frame $(X_1,X_2,X_3,X_4,X_5)$ will be called {\it 
strongly adapted} to $D$. 

We are going to show how to calculate our invariants 
starting from some adapted frame to distribution. Let again 
$\vec h=\vec h_{X_1,X_2}$ as in (\ref{ham25}). For any 
$\lambda\in (D^2)^\perp\backslash (D^3)^\perp$ consider the 
characteristic curve $\gamma$ of $D$ passing through 
$\lambda$. Under identification (\ref{tanid}) the reduced 
Jacobi curve $\widetilde J_\gamma$ lives in Lagrange 
Grassmannian $L(W_\lambda)$ of symplectic space 
$W_\lambda$, defined by (\ref{wl}). Let 
$\epsilon_1(\lambda)$ be the first vector in the canonical 
basis of the curve $t\mapsto \widetilde J_\gamma(e^{t \vec 
h}\lambda)$ at the point $t=0$. Note that it is more 
convenient to work directly with vector fields of 
$(D^2)^\perp$, keeping in mind that the symplectic space 
$W_\lambda$ belongs to the factor space 
$T_\lambda\bigl((D^2)^\perp\bigr)/ {\rm span} \bigl(\vec 
h(\lambda),\vec e(\lambda)\bigr)$. So, in the sequel by 
$\epsilon_1(\lambda)$ we will mean both the element of 
$W_\lambda$ and some representative of this element in 
$T_\lambda\bigl((D^2)^\perp\bigr)$, depending smoothly on 
$\lambda$. In the last case all equalities, containing 
$\epsilon_1(\lambda)$, will be assumed modulo ${\rm span} 
\bigl(\vec h(\lambda),\vec e(\lambda)\bigr)$. Now we are 
ready to prove the following 
\begin{prop}
\label{el25}
The vector $\epsilon_1(\lambda)$
can be chosen in the form
\begin{equation}
\label{e125eq} \epsilon_1(\lambda)= 6\Bigl(
\gamma_4(\lambda) \partial_{u_4}+\gamma_5(\lambda) 
\partial_{u_5}\Bigr), 
\end{equation}
where 
\begin{equation}
\label{gammarel} \gamma_4 (\lambda) u_5-\gamma_5(\lambda) 
u_4\equiv 1. 
\end{equation}
\end{prop}

{\bf Proof.}
First note that by (\ref{j1}) one has
\begin{equation}
\label{skewj1}
{\rm span}(\partial_{u_4},\partial_{u_5})\subset( {\cal J}^{1})^\angle.
\end{equation}
Hence by (\ref{dire1})
\begin{equation}
\epsilon_1=6\Bigl((\gamma_4\partial_{u_4}+\gamma_5\partial_{u_5}\Bigr)\,\,
\Bigl 
({\rm mod}\,{\rm span}(\vec h, \vec e)\Bigr), 
\end{equation}
where $\gamma_4 u_5-\gamma_5 u_4\neq 0$. Further, denote by 
$e_1(t)$ the first vector in the canonical basis of the 
curve $t\mapsto \widetilde J_\gamma(e^{t\vec h}\lambda)$. 
Then 
\begin{equation}
\label{cane1J}
e_1(t)=(e^{-t\vec h})_*\epsilon\bigl(e^{t\vec h}(\lambda)\bigr).
\end{equation}
Hence by (\ref{addt})
\begin{equation}
\label{barss}
\bar\sigma(e_1'(t),e_1''(t))\Bigr|_{t=0}\Bigl.=
\sigma\Bigl(\bigl[\vec h,[\vec h,\epsilon_1]\bigr](\lambda),[\vec h,\epsilon_1]
(\lambda)\Bigr).
\end{equation}
By direct computation one can show that 
\begin{eqnarray}
&~&[\vec h,\epsilon_1]=6\Bigl(\gamma_5 \vec u_1-\gamma_4\vec u_2+
(\gamma_4 u_4-\gamma_5 u_4)\partial_{u_3}\Bigr)\quad \Bigr({\rm mod}\,{\rm span}(\vec h, \vec e,\epsilon_1)\Bigr),\nonumber \\
&~& \label{she1}\\
&~&\bigl[\vec h [\vec h,\epsilon_1]\bigr]=6(\gamma_4 u_4-\gamma_5 u_4)(\vec u_3+u_4\partial_{u_1}+u_5\partial_{u_2}) \quad\Bigr({\rm mod}\,{\rm span}(\vec h, \vec e,\epsilon_1,[\vec h,\epsilon_1])\Bigr).\nonumber
\end{eqnarray}
From (\ref{she1}) it is easy to show that   
the righthand side 
of (\ref{barss}) is equal to $36 (\gamma_4 u_5-\gamma_5 
u_4)^2$, which together with (\ref{seprime}) implies 
(\ref{gammarel}). $\Box$ 
\medskip

As a direct consequence of the previous Proposition, 
Proposition \ref{eqe1}, and relations (\ref{addt}), 
(\ref{cane1J}), (\ref{projcurv1}) we obtain 
\begin{theor}
\label{fundway1}
Let $\epsilon_1(\lambda)$ be
as in (\ref{e125eq}) and (\ref{gammarel}). Then there exist functions $A_0$ ,$A_1$ on
$(D^2)^\perp$ such that
\begin{equation}
\label{ad4h} ({\rm ad} \vec h)^4(\epsilon_1)=A_0 \epsilon_1+ \vec
h(A_1) {\rm ad}\vec h(\epsilon_1)+ A_1({\rm ad} \vec
h)^2(\epsilon_1)\,\, {\rm mod}\Bigl({\rm span} \bigl(\vec h,\vec
e\bigr)\Bigr)
\end{equation}
The fundamental form $A(\lambda)$ and the projective Ricci curvature 
$\rho_D(\lambda)$ of the distribution $D$ satisfy:
\begin{equation}
\label{fund25eq} 35 A=A_0+\frac{9}{100} A_1^2-\frac{3}{10}(\vec
h)^2(A_1),
\end{equation}
\begin{equation}
\label{rho25}
 \rho{_{_D}}=\frac{-\frac {2}{15} A_1
A^2-\frac{1}{6}\vec h \circ\vec h (A) A+\frac{3 
}{16}\bigl(\vec h(A)\bigr)^2}{|A|^{\frac {5}{2}}}.
\end{equation}
\end{theor}

\begin{remark}
\label{mode1}
{\rm It is clear that in the previous theorem we can take 
$\epsilon_1$ satisfying (\ref{e125eq}) and the following 
relation 
\begin{equation}
\label{gammarel1} 
\gamma_4(\lambda) u_5-\gamma_5(\lambda) u_4\equiv const
\end{equation}
along any characteristic curve of $D$ (instead of (\ref{gammarel})).
In particular one can take as $\epsilon_1$ one of the 
following vector fields: 
\begin{equation}
\label{samples}
\frac{1}{u_5}\partial_{u_4},\,\,\, \frac{1}{u_4}\partial_{u_5},\,\,\, 
\frac{(u_5\partial_{u_4}-u_4\partial_{u_5})}{u_4^2+u_5^2},\,\,\, 
\frac{(u_5\partial_{u_4}+u_4\partial_{u_5})}{u_5^2-u_4^2}
\end{equation}} $\Box$
\end{remark}

The formulas (\ref{ad4h}),(\ref{fund25eq}), and (\ref{rho25}) give an 
explicit way to calculate the fundamental form  and projective Ricci curvature 
of distribution, starting from some adapted frame of the distribution.
We will demonstrate later these formulas on several examples.
The previous theorem allows to prove also the following theorem about the algebraic structure of $(2,5)$-distributions
\begin{theor}
\label{algstruct} 
 In the case $n=5$ for any $q\in M$ the
 restriction of the densities $A(\cdot)$ of the fundamental 
form to the fibers $(D^2)^\perp(q)$ are degree $4$ 
homogeneous polynomials on $(D^2)^\perp(q)$. 
\end{theor}

{\bf Proof.}
Let $\epsilon_1=\frac{1}{u_5}\partial_{u_4}$. Also denote
$$
\begin{array}{l}
\tilde {\cal X}= \vec u_2-u_5\partial_{u_3}\\
\tilde Y_j=\vec u_j+\sum_{i=1}^3\{u_i,u_j\}\partial_{u_i},\quad j=4,5
\end{array}$$
and let $F$ be as in (\ref{F25}). Then the tuple of vector fields
\begin{equation}
\label{frpol}
\epsilon_1,\tilde {\cal X}, F, Y_4,Y_5,\vec h, \vec e
\end{equation}
constitute a frame on $(D^2)^\perp$.
By direct calculations
\begin{eqnarray}
&~& [\vec h,\epsilon_1]=-\frac{1}{u_5}\tilde {\cal X}+ p_1 
\epsilon_1 \quad {\rm mod}\, \mathbb R \vec 
e\label{he1po}\\ &~&({\rm ad}\vec h)^2(\epsilon_1)= F+p_2 X 
+ p_3 \epsilon_1\quad {\rm mod}\Bigl({\rm span}(\vec h, 
\vec e)\Bigr),\label{hhe1po} 
\end{eqnarray}
 where  $p_i$, $i=1,2,3$, are some rational functions in $u_4,u_5$ with  
denominator of the form $u_5^l$. From the form of vector 
fields $\vec h$ and $e_1$ it follows that the coordinates 
of vector field $({\rm ad}\vec h)^4(\epsilon_1)$ w.r.t. the 
frame (\ref{frpol}) are also rational functions in 
$u_4,u_5$ with 
 denominator of the form $u_5^l$. But from (\ref{ad4h}), (\ref{he1po}), and
(\ref{hhe1po}) it follows that
$$({\rm ad}\vec h)^4\epsilon_1\subset {\rm span}(\epsilon_1,\tilde {\cal X},
F,\vec h,\vec e).$$
Expressing $\tilde{\cal X}$ and $F$ by $\epsilon_1$, $[\vec h,\epsilon_1]$, 
and 
$({\rm ad}\vec h)^2(\epsilon_1)$  from (\ref{he1po}) and (\ref{hhe1po}) ${\rm mod}({\rm span}(\vec h, \vec e))$, one obtains that coefficients $A_0$, $A_1$ from (\ref{ad4h}) and 
so also the fundamental form $A$ are rational functions in $u_4,u_5$ with  
denominator of the form $u_5^l$. But by Remark \ref{remsmooth} $A$ is smooth
at the points with $u_5=0$, $u_4\neq 0$. It implies that $A$ has to be 
polynomial. (Another argument is as follows:
if at the beginning one takes $\epsilon_1=\frac{1}{u_4}\partial_{u_5}$, 
then similarly to above one obtains that  $A$ is a rational function in $u_4,u_5$ with  denominator of the form $u_4^l$, which implies that $A$ has to be polynomial). $\Box$
\medskip
 
As a direct consequence of the previous theorem we obtain 

\begin{cor}
\label{tanpol} For any $q\in M$ the tangential fundamental 
form $\AA_q$ is degree $4$ homogeneous polynomial on 
$D(q)$. 
\end{cor}
\medskip

\begin{remark}
\label{compcar}
{\rm From the previous corollary it follows that
 the tangential fundamental form has the same algebraic nature, 
 as the {\it covariant binary biquadratic form}, constructed 
 by E.Cartan in \cite{cartan} (chapter VI, paragraph 33). 
 We call it Cartan's tensor.  In the next paper \cite{zelcart}
  we prove that the tangential fundamental form coincides   
(up to constant factor $-35$) with Cartan's tensor.} $\Box$ 
\end{remark}

\begin{remark}
{\rm In terms of canonical projective structure on abnormal 
extremal (see Remark \ref{bryantps}) and fundamental form 
one can obtain sufficient conditions for rigidity of the 
corresponding abnormal trajectory of $(2,5)$-distribution: 
A smooth curve $\xi$ tangent to distribution $D$ and 
connecting two points $q_0$ and $q_1$ is called rigid, if 
in some $C^1$-neighborhood of $\xi$ the only curves tangent 
to $D$ and connecting $q_0$ with $q_1$ are 
reparametrizations of $\xi$. Rigid curves are automatically 
abnormal trajectories of $D$. From result formulated in 
Remark \ref{rigid} (see \cite{zel}, Theorem 4.2 for precise 
statement) and comparison theorems from \cite{jac2} ( 
Theorem 5, item 1 there) one has the following:} For 
abnormal trajectory $\xi$ of $(2,5)$-distribution to be 
rigid it is sufficient the existence of a global projective 
parameter on $\xi$ together with the nonpositivity of the 
fundamental form along $\xi$ (equivalently nonnegativity of 
Cartan's tensor along $\xi$). 
 {\rm Moreover, if some Riemannian metric is given on $M$, then under the same conditions the corresponding abnormal trajectory is the shortest  among all curves tangent to distribution $D$, connecting 
its endpoints and sufficiently closed to this abnormal 
trajectory in $C$-topology. It follows again from the 
mentioned comparison theorem and from the fact that 
simplicity of the Jacobi curve of the abnormal extremal 
implies minimality of the length of the corresponding 
abnormal trajectory in $C$-topology (see \cite{agrmin} and 
\cite{agrc0}).}$\Box$ 
\end{remark}

{\bf 3.3 Examples.} Now we will give five examples of concrete distributions or families of distributions, for which we
 have computed the fundamental form and projective Ricci curvature using 
 Theorem \ref{fundway1}. We will present the computations 
only in Example 4, while in other case we will give only 
the results (in fact Examples 2 and 3 are included in 
Example 4; Examples 1-3 and other examples with the 
detailed computations can be found in \cite{zeldoc}). 

But before let us
introduce some notations. Let $X_1,X_2,X_3,X_4,X_5$ be an adapted frame to the distribution and $u_i$, $1\leq i\leq 5$, be
the corresponding quasi-impulses, defined by (\ref{quasi25}).
Suppose that this frame satisfies the following commutative
relations
\begin{equation}
\label{comrel25} [X_i,X_j]=\sum_{i=1}^n c_{ji}^k X_k
\end{equation}
Then the Hamiltonian vector fields $\vec u_i$, corresponding to
the functions $u_i$ satisfy
\begin{equation}
\label{vecui} \vec u_i=X_i+\sum_{j=1}^5\sum_{k=1}^5 c_{ji}^k u_k
\partial_{u_j}
\end{equation}
Therefore, the restriction of the vector field $\vec h$ on
\begin{equation}
\label{D2u} \bigl(D^2\bigr)^\perp\cap T^* M= \{\lambda\in
T^*M:u_1(\lambda)=u_2(\lambda)=u_3(\lambda)=0\}
\end{equation}
satisfies
\begin{equation}
\label{vechc25}
\begin{array}{c}
 \vec h=u_4\vec u_2-u_5\vec u_1= u_4 X_2-u_5 X_1+
\Bigl(c_{42}^4 u_4^2+(c_{42}^5-c_{41}^4)u_4 u_5-c_{41}^5
u_5^2\Bigr)\partial_{u_4}+\\ +\Bigl(c_{52}^4
u_4^2+(c_{52}^5-c_{51}^4)u_4 u_5-c_{51}^5
u_5^2\Bigr)\partial_{u_5}.
\end{array}
\end{equation}

{\bf Example 1. Free nilpotent (2,5)- distribution.} Let 
$L_1$ be the $5$-dimensional nilpotent Lie algebra with the 
following commutation rules in some basis $X_1,\ldots,X_5$: 
\begin{eqnarray}
&~& [X_1,X_2]=X_3,\,\,[X_1,X_3]=X_4,\,\, 
[X_2,X_3]=X_5\nonumber\\ &~& ad X_4=0,\,ad 
X_5=0\label{compnilp} 
\end{eqnarray}
Actually $L_1$ is the free nilpotent $3$-step Lie algebra 
with two generators. Let $M_1$ be 
the Lie group with the Lie algebra $L_1$. We consider 
$X_1,\ldots, X_5$ as left-invariant vector fields on $M_1$. 
Let $D_1=span(X_1,X_2)$. Such distribution is called {\it 
free nilpotent $(2,5)$-distribution}. 

By (\ref{compnilp}) the tuple of left invariants fields $(X_1,X_2,X_3,X_4,X_5)$ constitutes a strong adapted frame to distribution $D_1$. Applying 
Theorem {\ref{fundway1} to this frame, it is easy to show that 
 the fundamental form $A_{D_1}$ of distribution of $D_1$ satisfies
\begin{equation}
\label{AD1}
A_{D_1}\equiv 0
\end{equation}

{\bf Example 2. Left-invariant rank 2 distribution on $SO(3)\times 
\mathbb {R}^2$ .} Denote by $M_2$ the Lie group $SO(3)\times\mathbb 
{R}^2$. Let $L_2= so(3)\oplus \mathbb {R}^2$ be Lie algebra 
corresponding to Lie group $M_2$ . Suppose that $E_{ij}$ is $3\times 
3$ matrix such that its $(i,j)$th entry is equal to $1$ and all other 
entries equal to $0$. Take the following basis $a_1,a_2,a_3$ in 
$so(3)$: 
\begin{equation}
\label{basso3} a_1=E_{12}-E_{21}, a_2=E_{13}-E_{31},
a_3=E_{32}-E_{23}
\end{equation}
Then $a_1,a_2,a_3$ satisfy the following commutative relations:
\begin{equation}
\label{comso3} [a_1,a_2]=a_3,\,\,[a_2,a_3]=a_1,\,\,[a_3,a_1]=a_2
\end{equation}
Let $b_1,b_2$ be some basis of $\mathbb{R}^2$. Denote
\begin{equation}
\label{vecD2}
X_1=(a_1,b_1),\,\,X_2=(a_2,b_2),
\end{equation}
\begin{equation}
\label{D2} D_2={\rm span}\,(X_1,X_2) 
\end{equation}
We consider $X_1$, $X_2$,
as left-invariant vector fields on Lie group $M_2$. 
Consequently, $D_2$ defined by (\ref{D2}) can be considered 
as left-invariant rank 2 distribution on $M_2$.

\begin{remark}
\label{so3un} {\rm It can be shown easily (see \cite{zeldoc}) 
that the distribution $D_2$ is unique, up to group 
automorphism of $M_2$, left-invariant completely nonholonomic 
rank 2 distribution 
on $M_2$ and its small growth vector is $(2,3,5)$.} 
\end{remark}

\begin{remark}
\label{so3rol}
{\rm Note that distribution $D_2$ appears, when one
studies the problem of rolling ball on the plane without slipping and twisting
 (see Example 4 below
and also \cite{jurd} for the details).}
\end{remark}

Completing the basis $(X_1, X_2)$, defined by (\ref{vecD2}), to the strong adapted frame to $D_2$, and 
applying Theorem \ref{fundway1} to this frame, one has easily that 
the fundamental form $A_{D_2}$ and 
projective Ricci curvature $\rho_{D_2}$ of $D_2$ satisfy 

\begin{equation}
\label{fundso3}
 A_{D_2}\sim (u_4^2+u_5^2)^2,
\end{equation}
\begin{equation}
\label{rhoD2} \rho{_{_{D_2}}}=\frac{4\sqrt{35}}{9}
\end{equation}
(here as in the sequel we use the sign $\sim$ to emphasize that the fundamental form at the point is defined up to multiplication on a positive constant).
\begin{concl} From (\ref{AD1}) and (\ref{fundso3}) it
follows that germs of  distributions $D_1$ and $D_2$  are not
equivalent
\end{concl}

\begin{remark}
\label{nilpcart} {\rm Actually, the rank 2 distribution on 
$5$-dimensional manifold has the identically zero 
fundamental form iff it is locally equivalent to the 
distribution $D_1$. It follows from the fact that our 
fundamental form coincides with Cartan's tensor ( see \cite{zelcart}) and the fact that Cartan's 
tensor of distribution is identically zero iff it is 
locally equivalent to the distribution $D_1$ ( see chapter 
VIII of \cite {cartan}).}$\Box$ 
\end{remark}

 {\bf Example 3.
Left-invariant rank 2 distributions on
$SL(2,\mathbb{R})\times\mathbb {R}^2$.}
Let $E_{ij}$ be $2\times 2$ matrix such that its $(i,j)$th entry is 
equal to $1$ and all other entries equal to $0$. Take the following 
basis $a_1,a_2,a_3$ in $sl(2,\mathbb{R})$: 
\begin{equation}
\label{bassl2} a_1=\frac{1}{2}(E_{11}-E_{22}),
a_2=\frac{1}{2}(E_{12}-E_{21}), a_3=\frac{1}{2}(E_{12}+E_{21})
\end{equation}
Then $a_1,a_2,a_3$ satisfy the following commutative relations:
\begin{equation}
\label{comsl2} [a_1,a_2]=a_3,\,\,[a_2,a_3]=a_1,\,\,[a_3,a_1]=-a_2,
\end{equation}
 Denote by
 \begin{equation}
 \label{phep}
 P_h={\rm span}(a_1,a_2),\,\,\,P_e={\rm span}(a_1,a_3)
\end{equation}
Note that the restriction of the Killing form on $P_h$ is 
indefinite nondegenerated and on $P_e$ is positive definite 
quadratic form. Let $b_1,b_2$ be some basis of 
$\mathbb{R}^2$. Suppose that 
\begin{equation}
 \label{Dhep}
 D_{3,h}={\rm span}\bigl((a_1, b_1),(a_2,b_2)\bigr),\,\,\,
D_{3,e}={\rm span}\bigl((a_1,b_1),(a_3,b_2)\bigr),
\end{equation}
We consider $D_{3,h}$ and $D_{3,e}$, as left-invariant rank 2
distribution on the Lie group $SL(2,\mathbb{R})\times
\mathbb{R}^2$.

\begin{remark}
\label{sl2two} It can be shown easily (see \cite{zeldoc}) 
that distributions $D_{3,h}$ and $D_{3,e}$ are the only two 
different left-invariant rank 2 completely nonholonomic 
distributions on $SL(2,\mathbb{R})\times \mathbb{R}^2$, up 
to Lie group automorphisms of $SL(2,\mathbb{R})\times 
\mathbb{R}^2$, and their small growth vector is $(2,3,5)$. 
\end{remark}

\begin{remark}
\label{sl2rol} {\rm Note that distribution $D_{3,e}$ 
appears, when one studies the problem of rolling hyperbolic 
plane on the Euclidean plane without slipping and twisting, 
(see geometric model in Example 4 below).} 
\end{remark}
Completing the bases, chosen in (\ref{Dhep}), to the strong adapted frames 
of $D_{3,h}$ and $D_{3,e}$  and 
applying Theorem \ref{fundway1} to these frames, one has
easily that  the fundamental 
form $A_{D_{3,h}}$ and 
the projective curvature $\rho_{_{D_{3,h}}}$  of $D_{3,h}$ satisfy
\begin{equation}
\label{fundsl2h}
 A_{D_{3,h}}\sim (u_4^2-u_5^2)^2,
\end{equation}
\begin{equation}
\label{rhoD3h} \rho_{_{D_{3,h}}}=
\left\{
\begin{array}{cc}
-\frac{4\sqrt{35}}{9}& u_4>u_5\\
\frac{4\sqrt{35}}{9}& u_4<u_5
\end{array}\right.,
\end{equation}
while the fundamental form $A_{D_{3,e}}$ and 
the projective curvature $\rho_{_{D_{3,e}}}$ of  $D_{3,e}$ satisfy
\begin{equation}
\label{fundsl2e}
 A_{D_{3,e}}\sim(u_4^2+u_5^2)^2,
\end{equation}
\begin{equation}
\label{rhoD3e} \rho{_{_{D_{3,e}}}}=-\frac{4\sqrt{35}}{9}
\end{equation}

\begin{concl}
From (\ref{fundso3}) and (\ref{fundsl2h}) it follows that 
the germs of distributions $D_{3,h}$ and $D_2$ are not 
equivalent; from (\ref{fundsl2h}) and (\ref{fundsl2e}) it 
follows that the germs of distributions $D_{3,h}$ and 
$D_{3,e}$ are not equivalent; finally, from (\ref{rhoD2}) 
and (\ref{rhoD3e}) it follows that the germs of 
distributions $D_2$ and $D_{3,e}$ are not equivalent. 
\end{concl}

\begin{remark}
\label{rho>A}{\rm
Note that for distributions $D_2$ and $D_{3,e}$
the information about their fundamental forms does not 
imply  
their local nonequivalence: fundamental forms in both cases 
are squares of sign definite quadratic forms. This is the 
case when the projective Ricci curvature helps to distinct 
the distributions.} 
\end{remark}

{\bf Example 4. Rolling of two spheres without slipping and 
twisting.} Rank $2$ distributions on $5$-dimensional 
manifold appear naturally when one studies the possible 
motions of two surfaces $S$ and $\widehat S$ in $\mathbb 
R^3$, which roll one on another without slipping and 
twisting. Here we follow the geometric model of this 
problem given in \cite{rolbod} (this model ignores the 
state constraints that correspond to the admissibility of 
contact of the bodies embedded in $\mathbb {R}^3$). The 
state space of the problem is the 5-dimensional manifold 
$$M_4=\{B:T_x S\mapsto T_{\hat x}\widehat S|\,B {\rm 
\,\,is\,\, an\,\, isometry}\}.$$ Let $B(t)\subset M_4$ be 
an admissible curve, corresponding to the motion of the 
rolling surfaces. Let $x(t)$ and $\hat x(t)$ be 
trajectories of the contact points in $S$ and $\widehat S$ 
respectively (so $B(t)$ can be considered as an isometry 
from $T_{x(t)} S$ to $T_{\hat x(t)}\hat S$).The condition 
of absence of slipping means that 
\begin{equation}
\label{slip} B(t) \dot x(t)=\dot{\hat x}(t), \end{equation} 
while the condition of absence of twisting can be written 
as follows 
\begin{equation} \label{twist} B(t)({\rm 
vector\,\,field\,\, parallel\,\, along}\,\, x(t))= {\rm 
vector\,\,field\,\, parallel\,\, along}\,\, \hat x(t). 
\end{equation}
From conditions (\ref{slip}) and (\ref{twist}) it follows 
that a curve $x(t)\in S$ determines completely the whole 
motion $B(t)\in M_4$ and the velocities of admissible 
motions define a $(2,5)$-distribution $D_{4,S,\hat S}$ on 
$M_4$. If $(v_1,v_2)$ and $(\hat v_1,\hat v_2)$ are some 
local orthonormal frames on $S$ and $\hat S$ respectively 
and $\beta$ is the angle of rotation from the frame $(B 
v_1(x), B v_2(x))$ to the frame $(\hat e_1(\hat x),\hat 
e_2(\hat x))$, then the points of $M_4$ are parametrized by 
$(x,\hat x,\beta)$ and one can choose a local basis of
distribution $D_{4,S,\hat S}$ as follows
\begin{eqnarray}
&~&X_1=v_1+\cos\beta\hat v_1+\sin\beta\hat v_2- 
(-\sigma_1+\hat\sigma_1\cos\beta+\hat\sigma_2\sin\beta)\partial_\beta,\nonumber\\ 
&~&~~ \label{rolbas0}\\ &~&X_2=v_2-\sin\beta\hat 
v_1+\cos\beta\hat v_2 +
(-\sigma_2-\hat\sigma_1\sin\beta+\hat\sigma_2\cos\beta 
\partial_\beta,\nonumber 
\end{eqnarray}
where $\sigma_i$,$\hat\sigma_i$ are structural functions of 
the frames: 
\begin{equation} 
\label{sig} [v_1,v_2]=\sigma_1 v_1+ \sigma_2 v_2,\quad[\hat 
v_1,\hat v_2]=\hat\sigma_1\hat v_1+\hat \sigma\hat v_2 
\end{equation} 
 
 Let us restrict ourselves to the case, when 
$S$ and $\widehat S$ are spheres of radiuses $r$ and $\hat 
r$ respectively. We will denote the corresponding 
$(2,5)$-distribution by $D_{4,r,\hat r}$. 

\begin{remark}
\label{trivrr} {\rm Obviously, distributions $D_{4,r,\hat 
r}$ with the same ratio $\frac{\hat r}{r}$ are equivalent 
and distributions $D_{4,r,\hat r}$ and $D_{4,\hat r,r}$ are 
equivalent too. The question is whether distributions 
$D_{4,r,\hat r}$ with different ratios $\frac{\hat 
r}{r}\geq 1$ are equivalent. }
\end{remark}

It turns out that the calculation of fundamental form and 
projective Ricci curvature of $D_{4,r,\hat r}$ gives the 
answer to the question in the previous remark. Taking the 
spherical coordinates $(\varphi,\psi)$ on the sphere $S$ , 
where $\varphi$ is the "altitude" (with values between $0$ 
and $\pi$) and $\psi$ is the "longitude", one can choose 
the following orthonormal frame on $S$ 
\begin{equation} \label{orth}v_1={1\over r 
}\partial_\varphi, \quad v_2={1\over r \sin 
\varphi}\partial_\psi. 
\end{equation}

In the same way take the spherical coordinates on $\widehat 
S$ and the orthonormal frame $(\hat v_1,\hat v_2)$ , 
defined by putting the sign $\hat{~}$ over $r$, $\varphi$, 
and $\psi$ in (\ref{orth}). Then the structural functions 
$\sigma_i$, $\widehat \sigma_i$ satisfy 
\begin{equation} \label{orthstruct} 
\sigma_1=0,\, \sigma_2=-{\cot\varphi\over r};\quad
 \widehat\sigma_1=0,\, 
\widehat\sigma_2=-{\cot \hat\varphi\over \hat r}.
\end{equation} 
Substituting (\ref{orthstruct}) into (\ref{rolbas}) one 
gets 
\begin{eqnarray}
&~&X_1=v_1+\cos\beta\hat v_1+\sin\beta\hat v_2- 
{\cot\hat\varphi\over\hat 
r}\sin\beta\partial_\beta,\nonumber\\ &~&~~ 
\label{rolbas}\\ &~&X_2=v_2-\sin\beta\hat v_1+\cos\beta\hat 
v_2\left({\cot\varphi\over r}- {\cot\hat\varphi\over\hat 
r}\cos\beta\right)\partial_\beta.\nonumber 
\end{eqnarray}
Let $X_3$, $X_4$, $X_5$ be as in (\ref{x345s}). Then from 
(\ref{rolbas}) one can obtain by direct computation 

\begin{eqnarray}
&~&X_3=-{\cot\varphi\over r} X_2+\left(\frac{1}{\hat 
r^2}-\frac{1}{r^2}\right)\partial_\beta\nonumber\\
&~&X_4=-{\cot\varphi\over r} X_3+{1\over 
r^2\sin^2\varphi}X_2- \left(\frac{1}{\hat 
r^2}-\frac{1}{r^2}\right)\left(X_2-v_2- {\cot\varphi\over 
r}\right)\partial_\beta\label{rolX345}\\ 
&~&X_5=\left(\frac{1}{\hat 
r^2}-\frac{1}{r^2}\right)(X_1-v_1)\nonumber 
\end{eqnarray}
 
Hence for $r\neq \hat r$ the tuple of the fields 
$(X_1,X_2,X_3, X_4, X_5)$ constitutes a strong adapted 
frame of the distribution $D_{4,r,\hat r}$. 

\begin{remark}
\label{r=r} {\rm It is clear that for $r=\hat r$ the 
distribution $D_{4,r,\hat r}$ is integrable.}$\Box$ 
\end{remark}
 
Let us calculate the fundamental form and the projective 
Ricci curvature in this frame. First, from (\ref{rolbas}) 
and (\ref{rolX345}) one can obtain by direct calculations 
that 
\begin{eqnarray}
&~& [X_1,X_4]=-\left(\frac{1}{\hat 
r^2}+\frac{2}{r^2\sin^2\varphi}\right) 
\frac{\cot\varphi}{r}X_2-\left(\frac{1}{\hat r^2}- 
\frac{2}{r^2\sin^2\varphi}\right)X_3- 
\frac{\cot\varphi}{r}X_4,\nonumber\\ 
&~&[X_2,X_5]=-\left(\frac{1}{\hat 
r^2}-\frac{1}{r^2\sin^2\varphi}\right) 
\frac{\cot\varphi}{r}X_2-\left(\frac{1}{\hat 
r^2}+\frac{\cot^2\varphi}{r^2}\right)X_3- 
\frac{\cot\varphi}{r}X_4,\label{rolcomm}\\ 
&~&[X_2,X_4]=[X_1,X_5]=0.\nonumber 
\end{eqnarray}
Then by (\ref{comrel25}) and (\ref{vechc25})
\begin{equation}
\label{rolh}
\vec h=u_4 X_2-u_5 X_1- u_4\frac{\cot\varphi}{r}\partial_\theta,
\end{equation}
where $\partial_\theta$ is as in (\ref{theta25}). Take the 
polar coordinates $u_4=R\cos\theta$, $u_5=R\sin\theta$ on 
the fibers $(D^2)^\perp(q)$. From (\ref{rolh}) it is clear 
that $\vec h(R)=0$. Hence by Remark \ref{mode1} one can 
take as $\epsilon_1$ in Theorem \ref{fundway1} the field 
$\partial_\theta$. Besides, by homogeneity of fundamental 
form it is sufficient to restrict our calculations to the 
set $\{R=1\}$. The vector field $\vec h$ has on this set 
the form 
\begin{equation}
\label{rolhpol}
\vec h=\cos \theta X_2-\sin\theta X_1- \cos\theta\frac{\cot\varphi}{r}\partial_\theta,
\end{equation}
Further, by direct computation, it is not hard to get
\begin{eqnarray}
&~&[\vec h,\partial_\theta]=\sin \theta X_2+\cos\theta X_1- 
\sin\theta\frac{\cot\varphi}{r}\partial_\theta,\nonumber\\
&~&({\rm ad}\vec h)^2(\partial_\theta)=-\left(X_3-\frac{\cot\varphi}{r}+\frac{1}{r^2}\partial_\theta\right),\label{adrol}\\
&~&({\rm ad}\vec h)^3(\partial_\theta)=-\cos\theta X_5+\sin\theta X_4+
\sin\theta\frac{\cot\varphi}{r} X_3-\frac{\sin\theta}{r^2\sin^2\varphi}X_2-\frac{1}{r^2}[\vec h,\partial_\theta].\nonumber
\end{eqnarray}
Finally, from (\ref{rolcomm}) and (\ref{adrol}) one can 
obtain without difficulties that 
\begin{equation}
\label{adrol4}
({\rm ad}\vec h)^4(\partial_\theta)=-\frac{1}{r^2\hat r^2}\partial_\theta-
\left(\frac{1}{r^2}+\frac{1}{\hat r^2}\right)
({\rm ad}\vec h)^2(\partial_\theta),
\end{equation}
which together with (\ref{fund25eq}) and (\ref{rho25}) 
implies that the fundamental for $A_{D_{4,r,\hat r}}$ and projective Ricci curvature $\rho_{_{D_{4,r,
\hat r}}}$ of $D_{4,r,\hat r}$ satisfy
\begin{eqnarray}
&~& A_{D_{4,r,\hat r}}\sim {\rm sgn}\, 
\Bigl((9\hat r^2-r^2)(\hat r^2-9 r^2)\Bigr)(u_4^2+u_5^2)^2,\label{Arol}\\
&~&\rho_{_{D_{4,r,\hat r}}}=\frac{4 \sqrt{35}}{3}\frac{r^2+\hat r^2}
{\sqrt{|(9\hat r^2-r^2)(\hat r^2-9 r^2)|}}\label{rhorol}.
\end{eqnarray}

\begin{prop}
\label{noeqrr} Distributions $D_{4,r,\hat r}$ with 
different ratios $\frac{\hat r}{r}\geq 1$ are not 
equivalent. 
\end{prop}

{\bf Proof.} Let $$I_1=(1,3),\quad I_2=\{3\},\quad 
I_3=(3,+\infty).$$ From (\ref{Arol}) distributions 
$D_{4,r,\hat r}$ with ratios, taking values in different 
$I_j$, are not equivalent. Note that the function 
$\alpha\mapsto \frac{\alpha^2+1} {\sqrt{|( \alpha^2-9)( 
9\alpha^2- 1)|}}$ is monotone on both intervals $I_1$ and 
$I_3$. This together with (\ref{rhorol}) implies that 
distributions $D_{4,r,\hat r}$ with ratios, taking values 
in one of $I_j$, $j=1$ or $3$, are not equivalent, which 
together with Remark \ref{r=r} concludes the proof of the 
Proposition. $\Box$ 

\begin{remark}
\label{r3r}{\rm From (\ref {Arol}) and Remark 
\ref{nilpcart} it follows that if $\frac{\hat r}{r}=3$ or 
$\frac{1}{3}$, then $D_{4,r,\hat r}$ is locally equivalent 
to the free nilpotent $(2,5)$-distribution $D_1$ from 
Example 1.}$\Box$ 
\end{remark}
\begin{remark}
\label{varrol} {\rm From (\ref{Arol}) and (\ref{rhorol}) it 
follows that for $\frac{\hat r}{r}>3$ the fundamental form 
is square of sign definite quadratic form and the 
projective Ricci curvature varies from 
$\frac{4\sqrt{35}}{9}$ to $\infty$ ($\frac{4\sqrt{35}}{9}$ 
corresponds to $\frac{\hat r}{r}=\infty$, i.e., to the 
rolling of the sphere on the plane) , while for 
$1<\frac{\hat r}{r}<3$ the fundamental form is $-$ square 
of sign definite quadratic form and the projective Ricci 
curvature varies from $\frac{\sqrt{35}}{3}$ to $\infty$ 
(not including $\frac{\sqrt{35}}{3}$).} $\Box$ 
\end{remark}

More generally, if $S$ and $\widehat S$ are surfaces of 
constant curvatures $k$ and $\hat k$ respectively, then the 
fundamental form and projective Ricci curvature of the 
distribution, generated by the rolling of one surface on 
another without slipping and twisting, satisfy 

\begin{eqnarray}
&~& A_{D_{4,r,\hat r}}\sim {\rm sgn}\, 
\Bigl((9\hat k-\hat k)(\hat k-9 k)\Bigr)(u_4^2+u_5^2)^2,\label{Arol1}\\
&~&\rho_{D}=\frac{4 \sqrt{35}}{3}\frac{k+\hat k}
{\sqrt{|(9\hat k-k)(\hat k-9 k)|}}\label{rhorol1}.
\end{eqnarray}

{\bf Example 5. Distributions generated by curves of 
constant torsion on $3$-dimensional manifold of constant 
curvature.} These distributions were mentioned already in 
\cite{cartan} (chapter XI, paragraphs 52, 53). Let $Q$ be 
an oriented $3$-dimensional Riemannian manifold. Then for 
given $\tau$ the curves of constant torsion $\tau$ together 
with their binormals are admissible curves of a rank $2$ 
distribution on $5$-dimensional manifold $M_5=Q\times S^2$. 
Indeed, let $\gamma(t)$ be the curve in $Q$ without 
inflection points, and let $n(t)\in S^2$ be the 
corresponding binormal. Then $\gamma$ has a constant 
torsion $\tau$ iff 
\begin{equation}
\label{constor} \dot \gamma(t)=\frac{1}{\tau} n(t)\times 
\nabla_{\dot\gamma(t)}n(t), \end{equation} where by 
$\times$ we mean the vector product induced on each (oriented) 
tangent space $T_{\gamma(t)}Q$ by the Riemannian metric and $\nabla$ 
denotes the covariant derivative, corresponding to this 
metric. Obviously, relation (\ref{constor}) defines the 
rank 2 distribution on $M_5$. We restrict ourselves to the 
case when $Q$ has constant curvature $K$ and denote by 
$D_{5,\tau,K}$ the corresponding $(2,5)$-distribution. It 
can be shown that the corresponding fundamental form 
$A_{_{D_{5,\tau, K}}}$ and the projective Ricci curvature 
$\rho_{_{D_{5,\tau, K}}}$ satisfy 
\begin{eqnarray} &~& \label{Akruch} A_{_{D_{5,\tau, 
K}}}\sim {\rm 
sgn}\Bigl((\frac{\tau^2}{K}-4)(1-4\frac{\tau^2}{K})\Bigr)\Bigr 
({\rm sign\,\,definite\,\,quadratic\,\,form}\Bigl)^2\\ 
&~&\label{rhokruch} \rho_{_{D_{5,\tau, 
K}}}=\frac{2\sqrt{35}}{3}\frac{\frac{\tau^2}{K}+1} 
{\sqrt{|(\frac{\tau^2}{K}-4)(1-4\frac{\tau^2}{K})|}} 
\end{eqnarray}
 
Suppose that $S$ is three-dimensional sphere of radius $R$. 
Then $K=\frac{1}{R^2}$. Note that the expressions in 
(\ref{Akruch}) and (\ref{rhokruch}) are invariant w.r.t. 
transformation $\tau\mapsto \frac{1}{R^2\tau}$. 
\begin{remark}
\label{nilpkruch} {\rm If $\tau R=2$ or $\frac{1}{2}$, then 
from (\ref{Akruch}) the fundamental form is equal to zero. 
Hence by Remark \ref{nilpcart} the corresponding 
distribution is locally equivalent to the free nilpotent 
$(2,5)$-distribution $D_1$ and by Remark \ref{r3r} it is 
locally equivalent to the distribution $D_{4,r,\hat r}$ 
with $\frac{\hat r}{r}=3$ or $\frac{1}{3}$.}$\Box$
\end{remark}

\begin{remark}
\label{sqkruch} {\rm If $0<\tau R 
<\frac{1}{2}$ or $\tau R>2$ then the fundamental form is 
$-$ square of sign definite quadratic form and the 
projective Ricci curvature varies monotonically on both 
intervals from $\frac{\sqrt{35}}{3}$ to $\infty$ ( 
$\frac{\sqrt{35}}{3}$ corresponds to the case, when $Q$ is 
Euclidean space). If $\frac{1}{2}<\tau R 
< 1$ or $1<\tau R<2$ then the fundamental 
form is square of sign definite quadratic form and the 
projective Ricci curvature varies monotonically on both 
intervals from $\frac{4 \sqrt{35}}{9}$ to $\infty$ ( not 
including $\frac{4 \sqrt{35}}{9}$).}$\Box$ \end{remark} 

Note that in the case $\tau=0$ the distribution is 
integrable, while in the case $\tau R=1$ the square of the 
distribution is rank 3 integrable. \medskip 

Till now we used our invariants in order to prove the 
nonequivalence of distributions. But what to do, if both 
the fundamental form and the projective Ricci curvature do 
not distinct distributions? For example by Remarks 
\ref{varrol} and \ref{sqkruch} for any ratio $\frac{\hat 
r}{r}$ there exists distribution 
$D_{5,\tau,\frac{1}{R^2}}$, which has the fundamental form 
of the same type and the same projective Ricci curvature as 
$D_{4,r,\hat r}$. Does it imply that these distributions 
are equivalent? We will treat the questions of this kind in 
the forthcoming paper \cite{zelprep}. Below we formulate a 
theorem, which will be proved in this paper:

\begin{theor}
\label{clastheor1} For given $s\in\{1,-1\}$ and $\rho\in R$ 
there exists unique, up to diffeomorphism, germ of 
 $(2,5)$-distribution satisfying the following three conditions: 

1. Its fundamental form is $s$ multiplied by the square of 
a nondegenerated quadratic form $Q$; 

2. Its symmetry group is $6$-dimensional; 

3. If $Q$ is sign definite, then its projective Ricci 
curvature is identically equal to $\rho$, if $Q$ is sign 
indefinite, then the absolute value of its projective Ricci 
curvature is identically equal to $|\rho|$. 
\end{theor}

\begin{remark}
\label{dim6} {\rm It can be shown that if distribution $D$ 
satisfies condition 1 of Theorems \ref{clastheor1}, 
then the dimension of the group of 
symmetries of $D$ is not greater than $6$.}$\Box$ 
\end{remark}

\begin{remark}
{\rm It can be shown that conditions 1 and 2 of Theorem 
\ref{clastheor1}  imply that 
projective Ricci curvature or its absolute value is 
identically equal to some constant.}$\Box$
\end{remark}

It is easy to see that the group of symmetries of 
distribution $D_{4,r,\hat r}$ contains a subgroup 
isomorphic to $SO(3)\times SO(3)$ and therefore by Remark 
\ref{dim6} it is $6$-dimensional for $\frac {\hat r}{r}\neq 
3$ or $\frac{1}{3}$, while the group of symmetries of 
distribution $D_{5,\tau,\frac{1}{R^2}}$ contains a subgroup 
isomorphic to $SO(4)$ and therefore by Remark \ref{dim6} it 
is also $6$-dimensional for $\tau R\neq 2$ or 
$\frac{1}{2}$. Therefore Theorem \ref{clastheor1} implies 
\begin{cor}
\label{D45} If distributions $D_{4,r,\hat r}$ and 
$D_{5,\tau,\frac{1}{R^2}}$ have the fundamental forms, 
which are $\pm$ square of positive definite quadratic form, 
where the sign is the same for both distributions, and 
their projective Ricci curvatures are the same constant, 
then these distributions are locally equivalent. 
\end{cor}

\section{Algebraic structure of fundamental form in the 
case $n>5$}
 \indent\setcounter{equation}{0}

In the present section we show that in the case $n>5$ the 
fundamental form is in general a rational function, which 
is not a polynomial. By Remark \ref{remsmooth} 
singularities of fundamental form could occur out of the 
set ${\cal R}_D$, i.e., at the points, where the weight of 
the corresponding Jacobi curve is not constant. 

First take some curve $\Lambda (t)$ in the
 Grassmannian of half-dimensional subspaces 
$G_m(W)$ such that at some point $\bar t$ it has a jump of 
the weight and look more carefully what happens with Ricci 
curvature and fundamental form of this curve near this 
point. More precisely, we consider the following situation: 
in some punctured neighborhood of $\bar t$ the curve 
$\Lambda(t)$ has constant weight $k$, while at point $\bar 
t$ it has the weight $k+1$. In this case we will say that 
$\bar t$ is {\it a point of the weight jump one} of 
$\Lambda(t)$. 

\begin{lemma}
\label{jump1lem} If the curve $\Lambda(t)$ in $G_m(W)$ has 
a point of the weight jump one at $\bar t$, then the Ricci 
curvature has a pole of order 2 at $t=\bar t$. If in 
addition the weight of $\Lambda(t)$ in the punctured 
neighborhood of $\bar t$ is greater than $1$, then the 
density of fundamental form of this curve has a pole of 
order $4$ at $t= \bar t$. 
\end{lemma} 

{\bf Proof.} Suppose that $k$ is the weight of $\Lambda(t)$ 
in the punctured neighborhood of $\bar t$. Then the 
following function 
\begin{equation}
\label{funX} 
X(t_0,t_1)=\frac{\det(S_{t_0}-S_{t_1})}{(t_0-t_1)^k} 
\end{equation}
is smooth. By assumptions about the jump of the weight at 
$\bar t$ one has 
\begin{equation}
\label{Xjump1}
 X(\bar t, \bar t)=0,\quad 
{\partial\over\partial t_0}X(t_0,t_1)\Bigl|_{t_0=t_1=\bar 
t}\Bigr.= {\partial\over\partial 
t_1}X(t_0,t_1)\Bigl|_{t_0=t_1=\bar t}\Bigr.\neq 0. 
\end{equation}
 It implies that the function $X(t_0,t_1)$ is symmetric. 
 Indeed by permuting $t_0$ and $t_1$ in (\ref{funX}) we 
obtain that $X$ can be either symmetric or antisymmetric, 
but the last case is impossible, because $X(t,t)$ is not 
identically zero.

Without loss of generality it can be assumed that $\bar 
t=0$. From (\ref{Xjump1}) and symmetricity it follows that 
\begin{equation}
\label{Xfact} X(t_0,t_1)=(t_0+t_1) a(t_0,t_1), 
\end{equation}
where $a(t_0,t_1)$ is smooth. Using (\ref{Xfact}) and 
(\ref{gdet}) we obtain that the generating function 
$g_\Lambda(t_0,t_1)$ of the curve $\Lambda(t)$ satisfies 
\begin{equation}
\label{gX} g_\Lambda(t_0,t_1)={\partial^2 \ln\, 
X(t_0,t_1)\over\partial t_0\partial t_1} = 
-\frac{1}{(t_0+t_1)^2} +O(1) 
\end{equation}
It yields first that the Ricci curvature $\rho 
(t)=g_\Lambda(t,t)$ has the following expansion at $t=0 $ 
\begin{equation}
\label{rhojump} \rho(t)=g_\Lambda(t,t)=-\frac{1}{4 t
^2}+O(1) 
\end{equation}
In order to find an asymptotic for the density $A(t)$ of 
fundamental form at $t=\bar t$, we use formula 
(\ref{densg}). Again from (\ref{gX}) it follows easily that 
$${\partial^2\over 
\partial t_0^2}g_\Lambda(t_0,t_1)\Bigl|_{t_0=t_1=t}\Bigr.=
-\frac{3}{8 t^4}+O(1)$$ Substituting this into 
(\ref{densg}) and using (\ref{rhojump}) one gets 
\begin{equation}
\label{Ajump} A(t)= \frac{3(k-1)}{80\,k}\frac{1}{t 
^4}+O\left(\frac{1}{t^2}\right),\end{equation} which 
completes the proof of the Lemma. $\Box$ 
\medskip

\begin{remark}
\label{jump>1} {\rm If at some point of the curve in 
$G_m(W)$ the weight has jump greater than $1$, then in 
general the Ricci curvature and the density of fundamental 
form have also singularities, but the coefficients of the 
principal negative power in their Laurent expansion are not 
universal, as in (\ref{rhojump}) and (\ref{Ajump}).} $\Box$ 
\end{remark}
 
 The following proposition gives a simple 
characterization of the points of the weight jump one of 
the rank 1 curve $\Lambda(t)$ in Lagrange Grassmannian in 
terms of the subspaces ${\cal D}^{(i)}\Lambda$ defined by 
(\ref{primeiext}): 

\begin{prop}
\label{jump1prop}

The point $\hat t$ is the point of the weight jump one of 
rank 1 curve $\Lambda(t)$ in Lagrange Grassmannian $L(W)$, 
${\rm dim} W=2m$, iff the following relations holds 
\begin{equation}
\label{reljump1} {\rm dim}\, {\cal D}^{(m-1)}\Lambda(\bar 
t)={\rm dim}\, {\cal D}^{(m)}\Lambda(\bar t)=2m-1,\quad 
{\rm dim}\, {\cal D}^{(m+1)}\Lambda(\bar t)=2m 
\end{equation}
\end{prop}

The proof of the proposition can be easily obtained by 
application of some formulas and statements of section 6 
and 7 of \cite{jac1} (for example, formulas (6.15), (6.16), 
(6.18), (6.19), Lemma 6.1 and Proposition 3 there).

 Let us apply  Lemma \ref{jump1lem} and Proposition \ref{jump1prop}
 to distributions. 
 For this let subspace ${\cal J}^{(i)}(\lambda)$ 
be as in (\ref{Ji}). Denote by $S^0_D$ the following subset 
of $(D^2)^\perp\backslash \Bigl((D^3)^\perp\cup {\cal R}_D 
\Bigr) $ 
\begin{equation}
\label{sing0} S^0_D=\left \{\lambda\in (D^2)^\perp 
\backslash(D^3)^\perp 
:
\begin{array}{c} 
{\rm dim}\, {\cal J}^{(n-4)}(\lambda)={\rm dim}\,{\cal 
J}^{(n-3)}(\lambda)=2n-5,\\{\rm dim}\, {\cal 
J}^{(n-3)}(\lambda)=2n-4 
\end{array}
\right \}. 
\end{equation}  
Then by Proposition \ref{jump1prop} the set $S^0_D$ 
coincides with the subset of 
$(D^2)^\perp\backslash(D^3)^\perp$, consisting of points, 
in which the corresponding reduced Jacobi curves have the 
weight jump one. Also, from Propositions \ref{JT} one has

\begin{prop}
\label{Tjump} 
 The reduced Jacobi curve of the regular abnormal extremal 
$\gamma$ has the weight jump one at point $\lambda$ iff 
\begin{equation} \label{Tjumpeq} 
  \dim {\cal T}_\xi^{(n-4)}(q) =  \dim {\cal T}_\xi^{(n-3)}(q)=n-2,\quad  \dim {\cal 
    T}_\xi^{(n-2)}(q)=n-1
 \end{equation}
 where $\xi=\pi(\gamma)$ is the abnormal trajectory corresponding to 
 $\gamma$ and $q=\pi(\lambda)$.
 \end{prop}

Further note that by Proposition \ref{constJprop} the 
weight of these curves in the punctured neighborhoods of 
these points is equal to $(n-3)^2$ and therefore it is 
greater than $1$ in the considered cases. Suppose also that 
\begin{equation}
\label{sing0q}
S^0_D(q)=S^0_D\cap T_q^*M, \quad q\in M
\end{equation}
As a direct consequence of Lemma \ref{jump1lem} and 
Proposition \ref{rathom4} we obtain the following 

\begin{prop}
\label{singtheor} If the sets ${\cal R}_D(q)$ and 
$S^0_D(q)$ are not empty, then the fundamental form of 
distribution $D$ at point $q$ is a rational 
function, which is not a polynomial: all points of $S^0_D(q)$ are the points of discontinuity of it. 
\end{prop}

{\bf Example.} Consider distribution $\widetilde D$ in 
$\mathbb {R}^6$ spanned by the following vector fields 
\begin{equation}
\label{examp6} X_1={\partial\over\partial x_1}, \quad X_2= 
{\partial\over\partial x_2}+ x_1{\partial\over\partial 
x_3}+ \frac{x_1^2}{2}{\partial\over\partial 
x_4}+\Bigl(\frac{x_1^4}{4!}+\frac{x_1^2 x_2}{2}\Bigr)
{\partial\over\partial x_5}+x_1x_2{\partial\over\partial 
x_6}. 
\end{equation}  
Distribution $\widetilde D$ has the maximal possible small 
growth vector $(2,3,5,6)$ at any point. We claim that its 
fundamental form at $0$ is a rational function, which is not a polynomial. 
Indeed, by Proposition \ref{RD6prop}, the set ${\cal 
R}_{\widetilde D}(0)$ is not empty. It is not hard to show 
that the curve $(x_1,0,\ldots,0)$ is regular abnormal 
trajectory of corank 1. Moreover, from Proposition 
\ref{Tjump} it follows that the reduced Jacobi curve of any 
its lift $\gamma$ has the weight jump one at the point of 
intersection of $\gamma$ with $(\widetilde D^2)^\perp(0)$, 
which implies that $S^0_{\widetilde D}(0)$ is not empty. 
Now our claim follows from Proposition \ref{singtheor}. 
$\Box$ \medskip 

 From Proposition \ref{5weiconst} it is clear that for 
$(2,5)$-distributions with small growth vector $(2,3,5)$ 
the set $S^0_D$ is empty. For $n>5$ the situation is 
different: it turns out that after an appropriate 
complexification, one can prove that $S^0_D(q)$ is not 
empty for generic germ of $(2,n)$-distribution at $q$. 
Below we briefly describe this process, leaving the details 
to the reader. Note that for sufficiently big natural $l$ 
the weight, the jump of the weight and the density of 
fundamental form of the curve in Lagrange Grassmannian at 
some point is completely determined by its $l$-jet at this 
point. So for our purposes it is sufficient to work with 
$l$-jets of the curves with sufficiently big $l$ instead of 
the curves themselves. From the proof of Proposition 
\ref{rathom4} it follows that if for given $q\in M$ we 
consider a mapping $\Phi_l$ from $(D^2)^\perp(q)$ to the 
space of $l$-jets of curves in Lagrange Grassmannian, which 
assigns to any $\lambda\in (D^2)^\perp(q)$ the $l$-jet of 
the corresponding reduced Jacobi curve at $\lambda$, then 
the mapping $\Phi_l$ is rational on $(D^2)^\perp(q)$. Hence 
this mapping can be rationally continued from the real 
linear space $(D^2)^\perp(q)$ to its complexification 
$(D^2)^\perp(q)^{\mathbb C}$. After this continuation for 
any $\lambda\in (D^2)^\perp(q)^{\mathbb C}$ one has the 
$l$- jet of curve in the complex Lagrange Grassmannian 
(i.e., the set of (complex) half-dimensional subspace in 
complex even-dimensional space provided with skew-symmetric 
nondegenerated bilinear form). The theory of such curves is 
completely the same as in the case of real Lagrange 
Grassmannian. In particular, one has result analogous to 
Lemma \ref{jump1lem}. If we denote by $S^0_D(q)^{\mathbb 
C}$ the set of all $\lambda\in(D^2)^\perp(q)^{\mathbb C}$ 
such that the $l$-jet $\Phi_l(\lambda)$ corresponds to 
curve in the complex Lagrange Grassmannian with the weight 
jump one at $\lambda$, then we have the following 
generalization of Proposition \ref{singtheor}: 

\begin{prop}
\label{singtheor1} If the sets ${\cal R}_D(q)$ and 
$S^0_D(q)^{\mathbb C}$ are not empty, then the fundamental 
form of distribution $D$ at point $q$ is a rational 
function, which is not a polynomial: all points of $S^0_D(q)^{\mathbb C}$ are singular 
points of the analytic continuation of the fundamental form 
to $(D^2)^\perp(q)^{\mathbb C}$. 
\end{prop}

The set $S^0_D(q)^{\mathbb C}$ can be described in more 
constructive way. For this note that the mappings 
$\lambda\mapsto J^{(i)}(\lambda)$, $\lambda\in 
(D^2)^\perp(q)$, depend rationally on $\lambda$ and 
therefore can be rationally continued to 
$(D^2)^\perp(q)^{\mathbb C}$ (after this continuation we 
look on $J^{(i)}(\lambda)$ as on complex linear spaces). 
Then from construction it follows that 

\begin{equation}
\label{sing0C} S^0_D(q)^{\mathbb C}=\left \{\lambda\in 
(D^2)^\perp(q)^{\mathbb C} 
:\begin{array}{c} 
{\rm dim} \,{\cal J}^{(n-4)}(\lambda)={\rm dim}\, {\cal 
J}^{(n-3)}(\lambda)=2n-5,\\ {\rm dim}\, {\cal 
J}^{(n-3)}(\lambda)=2n-4 
\end{array}
\right\}
\end{equation}  
(here all dimensions are over $\mathbb C$).
\begin{prop}
\label{S0nonempty}In the case $n>5$ for generic germ of 
$(2,n)$-distribution $D$ at $q_0$ the set 
$S^0_D(q_0)^{\mathbb C}$ is not empty. 
\end{prop}

{\bf Proof.} Choose again a local frame $\{X_i\}_{i=1}^n$ 
on $M$ such that $X_1,X_2$ constitute a local basis of $D$ 
and $X_3$, $X_4$, $X_5$ satisfy (\ref{x345}). From 
Proposition \ref{genconct} the set 
$(D^2)^\perp(q_0)\backslash {\cal R}_D(q_0)$ can be 
represented as a zero level set of some polynomial in 
$u_j$, $j =4,\ldots,n$. Denote this polynomial by ${\cal 
P}_D$. Using (\ref{singset}) and definition of subspaces 
$J^{(i)}$ one can show without difficulties that the 
polynomial ${\cal P}_D$ is either homogeneous polynomial of 
degree $d_n$ or identically equal to zero. In the case 
$n>5$ the degree $d_n>0$, while $d_5=0$. For example, in 
the case $n=6$ from Proposition \ref{RD6prop} it follows 
that one can take ${\cal P}_D=\alpha_6$, where $\alpha_6$ 
is as in (\ref{a6mal}). Hence $d_6=2$. In general 
$d_n=\frac{(n-4)(n-3)}{2}-1$. 

 From definition of ${\cal P}_D$ and  relations (\ref{singset}),
(\ref{sing0C}) it is easy to see that 
\begin{equation}
\label{HP} S^0_D(q_0)^{\mathbb C}=\left \{\lambda\in 
(D^2)^\perp(q_0)^{\mathbb C} :{\cal P}_D(\lambda)=0,\,\vec 
h_{X_1,X_2} ({\cal P}_D)(\lambda)\neq 0\right\}, 
\end{equation}
where $\vec h_{X_1,X_2}$ is as in (\ref{ham25}). In other 
words, the set $S^0_D(q_0)^{\mathbb C}$ is empty iff the 
following condition holds $$ {\cal 
P}_D(\lambda)=0\,\Rightarrow\vec h_{X_1,X_2} ({\cal 
P}_D)(\lambda)= 0.$$ From the form of the vector field 
$\vec h_{X_1,X_2}$ it follows that the polynomial 
$h_{X_1,X_2} ({\cal P}_D)(\lambda)$ is homogeneous of 
degree $d_n+1$. 

In general, if we consider the space of pairs of 
polynomials $\bigl(p_1(\lambda),p_2(\lambda)\bigr)$, where 
$p_i$ are polynomials of the fixed degrees $s_i$,
then the set of all pairs 
$\bigl(p_1(\lambda),p_2(\lambda)\bigr)$ such that 
$p_1(\lambda)=0$ implies $p_2(\lambda)=0$ is a finite union 
of algebraic varieties of the whole space of pairs. From 
this and the fact that coefficients of the polynomials 
${\cal P}_D(\lambda)$, $\vec h_{X_1,X_2} ({\cal 
P}_D)(\lambda)$ are polynomials in some jets space of 
$(2,n)$- distributions, we conclude that 
 there 
exists an open set $\widetilde {\cal U}_n$ in Zariski 
topology of this jets space 
such that the set $S^0_D(q_0)$ is not empty iff the 
corresponding jet of $D$ at $q_0$ belongs to $\widetilde 
{\cal U}_n$. 
 Therefore, in the same way as in Proposition 
 \ref{genconct},
  in order to prove our 
proposition it is sufficient to give an example of germ of 
$(2,n)$-distribution such that $S^0_D(q_0)$ is nonempty. As 
such example one can take distribution $\bar D$ spanned by 
the following vector fields 
\begin{equation} 
\label{examprd2} X_1={\partial\over\partial x_1}, \quad 
X_2= {\partial\over\partial x_2}+\sum_{i=1}^{n-4} 
\frac{x_1^i}{i!}{\partial\over\partial 
x_{i+2}}+\frac{x_1^{n-2}}{(n-2)!}{\partial\over\partial 
x_{n-1}}+x_1x_2{\partial\over\partial x_n}, 
\end{equation} 
where $(x_1,\ldots,x_n)$ are some local coordinates on $M$, 
$q_0=(0,\ldots,0)$. It is easy to see that the curve 
$(x_1,0,\ldots,0)$ is regular abnormal trajectory of corank 
1. Moreover, from Proposition \ref{Tjump} it follows that 
the reduced Jacobi curve of any its lift $\gamma$ has the 
weight jump one at the point of intersection of $\gamma$ 
with $(\bar D^2)^\perp(q_0)$, which implies that $S^0_{\bar 
D}(q_0)$ is not empty. This completes the proof. $\Box$ 
\medskip

  Finally, as a direct consequence of Propositions 
  \ref{genconct}, 
\ref{singtheor1} and \ref{S0nonempty} one has the following
 
\begin{theor}
\label{gensing} In the case $n>5$ a generic germ of 
$(2,n)$-distribution $D$ at $q$ has the fundamental form, 
which is a nonpolynomial rational function on $(D^2)^\perp(q)$. 
\end{theor}

\end{document}